\documentclass[a4paper]{amsart} % ,9pt
\usepackage{amsmath}
\usepackage{amssymb}%This gives us symbols like \square
\usepackage{verbatim}%This gives us the \begin{comment} environment.
\usepackage{amsthm}%This gives us Definition style of theorem.
\usepackage{mathrsfs}%This gives us \mathscr
\usepackage{graphicx}
\usepackage{mathtools}
\usepackage[colorlinks,pdfdisplaydoctitle,linkcolor=blue,citecolor=blue]{hyperref}
\usepackage{caption}
\usepackage{subcaption}
\usepackage{url}
\usepackage{epstopdf}
\usepackage{pgf,tikz}
\usepackage{todonotes}
\usepackage{bbm}
\usepackage[notref,notcite,color]{}
\usepackage{extarrows}
\usetikzlibrary{arrows}
\usepackage{algorithm,algpseudocode} 
\usepackage{multirow}

\usepackage{amssymb,amsthm}    % amssymb package contains more mathematical symbols
\usepackage{graphicx}          % graphicx package enables you to paste in graphics
\usepackage{amsmath}         % best maths package on offer
\usepackage{caption}
\usepackage{subcaption}
\usepackage{dsfont}
\usepackage{listings}
\usepackage{bbm}
\usepackage{color}
\usepackage{ dsfont }
\usepackage{enumerate}
 \usepackage{relsize}
\usepackage{xfrac}

\renewcommand{\Delta}{\triangle}

\definecolor{sw}{rgb}{0,0,0.7}

\definecolor{nc}{rgb}{0.7,0.1,0.45}

\definecolor{xt}{rgb}{0.7,0,0}

\definecolor{cls}{rgb}{0.1,0.45,0.1}

\def \Ee[#1]{\mathcal{E}^{\text{{#1}}}}
\def\R{\mathbb{R}}   
\def\N{\mathbb{N}}

\def\pa[#1,#2]{\frac{\partial {#1}}{\partial {#2}} }

\def\idom[#1,#2,#3]{\int_{#1}\hspace{1pt} {#2} \hspace{1pt} \text{d}{#3}}
\def\res[#1,#2]{\left.{#1}\right|_{#2}}

\def\var[#1,#2]{\langle \delta \mathcal{E}^{\text{{#1}}}({#2}),v\rangle}
\def\vars[#1,#2,#3]{\langle \delta^2\mathcal{E}^{\text{{#1}}}({#2})v,{#3}\rangle}
\def\vard[#1,#2,#3,#4]{\langle \delta\mathcal{E}^{\text{{#1}}}({#2})-\delta\mathcal{E}^{\text{{#3}}}({#4}),v\rangle}

\def\E{\mathbb{E}}

\def\N{\mathbb{N}}

\newcommand{\be}{\begin{equation}}
\newcommand{\en}{\end{equation}}
\newcommand{\ben}{\begin{equation*}}
\newcommand{\enn}{\end{equation*}}
\newcommand{\bea}{\begin{aligned}}
\newcommand{\ena}{\end{aligned}}

\def\ba#1\ena{\begin{align}#1\end{align}}
\def\ban#1\enan{\begin{align*}#1\end{align*}}

\theoremstyle{plain}
\newtheorem{thm}{Theorem}[section]

\newtheorem{lem}[thm]{Lemma}

\newtheorem{cor}[thm]{Corollary}

\newtheorem{assumption}[thm]{Assumption}

\newtheorem{remark}[thm]{Remark}
\newtheorem{example}[thm]{Example}

\numberwithin{equation}{section}

\begin{document}
%\title[PDE-constrained optimization with correlated data]{PDE-constrained optimization under uncertainty with correlated data}
%\title[Statistical learning in time-dependent linear systems]{Statistical learning in time-dependent linear systems}
%\title[Markov-chain Ensemble Kalman Inversion]{Markov-chain Ensemble Kalman Inversion}
\title[The Ensemble Kalman filter for dynamic inverse problems]{The Ensemble Kalman filter for dynamic \\ inverse problems}

\author[S. Weissmann] {Simon Weissmann}
\address{Institute of Mathematics, University of Mannheim, 68131 Mannheim, Germany}
\email{simon.weissmann@uni-mannheim.de}

\author[N. K. Chada] {Neil K. Chada}
\address{Department of Mathematics , City University of Hong Kong, 83 Tat Chee Ave, Hong Kong}
\email{neilchada123@gmail.com}

\author[X. T. Tong] {Xin T. Tong}
\address{Department of Mathematics, National University of Singapore, 119077, Singapore}
\email{mattxin@nus.edu.sg}

%\author[C. Schillings] {Claudia Schillings}
%\address{Mannheim School of Computer Science and Mathematics, University of Mannheim, 68131 Mannheim, Germany}
%\email{c.schillings@uni-mannheim.de}


\begin{abstract}
%In inverse problems, one is interested in estimating unknown model parameters from noisy observational data.
%Traditionally inverse problems are solved in the setting where one assumes a fixed forward operator describing the observation model. In this
%article, we consider the extension of this to settings where we have a %time-dependent
%dynamic forward model, which is motivated through applications
%in scientific computation and engineering. We consider this for a particular derivative-free optimizer, which is the ensemble Kalman inversion (EKI). We introduce
%and motivate a new methodology entitled dynamic-EKI, which is a particle-based method with a changing forward operator. %In the setting we now
%As result the observed data is correlated.
% We analyze our new method, where we present results related to the control of our particle system through its covariance structure. 
%This includes moment bounds and an ensemble collapse. This is required in order to demonstrate a convergence result, where we establish convergence in expectation. We test our theoretical findings for our dynamic-EKI on 2D Darcy flow partial differential equation.

In inverse problems, the goal is to estimate unknown model parameters from noisy observational data. Traditionally, inverse problems are solved under the assumption of a fixed forward operator describing the observation model. In this article, we consider the extension of this approach to situations where we have a dynamic forward model, motivated by applications in scientific computation and engineering. We specifically consider this extension for a derivative-free optimizer, the ensemble Kalman inversion (EKI). We introduce and justify a new methodology called \textcolor{black}{dynamic-EKI (DEKI)}, which is a particle-based method with a changing forward operator.
We analyze our new method, presenting results related to the control of our particle system through its covariance structure. This analysis includes moment bounds and an ensemble collapse, which are essential for demonstrating a convergence result. We establish convergence in expectation and validate our theoretical findings through experiments with \textcolor{black}{DEKI} applied to a 2D Darcy flow partial differential equation.

%\red{[Seems different from our real setups?]}

%Various numerics are provided to verify our theory for time-dependent EKI, which are a PDE-constrained optimization problem and a toy matrix problem.
\end{abstract}

\maketitle

\section{Introduction}\label{sec:intro}
The focus of this work is on the research area of inverse problems \cite{BB18,AMS10,AT87}, which is involves learning parameters, of quantities of interest from 
observations which are corrupted by noise. In many instances of inverse problems, one usually adopts a setting where the dynamics generating
the observations, i.e., the model of interest, is fixed and independent of time. However, there are highly relevant applications where the dynamics related to the forward operator change within each time frame, resulting in new observations. Such applications include geophysical sciences, numerical weather prediction
such as the \textcolor{black}{Navier-Stokes equation} and electrical resistivity tomography from thermodynamics \cite{MW06,SCI92}.  This motivates the use of a time-dependent forward operator within inverse problems, where the observation of the unknown of interest changes over time. Modifying traditional inverse problems in the setup
we described poses  computational and mathematical challenges, which has resulted in very limited literature. Our aim is to overcome and tackle these challenges, in the context
of inverse problems and, in particular, where we exploit a particular inverse problem methodology known as the ensemble Kalman inversion \cite{CIRS18,ILS13,LR09}. 
%\textcolor{blue}{To include references here}
\subsection{Preliminaries}
Before we introduce the notion of ensemble Kalman inversion (EKI), we present the mathematical formulation of an inverse problem. Throughout this manuscript we consider an underlying probability space $(\Omega,\mathcal F,\mathbb P)$. Given a set of noisy observations $u \in \R^p$, we are interested in recovering
some unknown parameter $z \in \R^d$, where the relationship between both is defined as
\begin{equation}
\label{eq:inverse_p}
u=S(z) + w, \quad w \sim \mathcal{N}(0,\Gamma)\,.
\end{equation}
Here, $S:\R^d\to \R^p$ is the forward operator and we assume our data is corrupted by additive Gaussian noise. 
We will use this assumption throughout this work. Commonly, inverse problems are ill-posed and require some regularization
scheme to produce numerical solutions. One method, which is the method of interest in this article, is the application of EKI \textcolor{black}{as a derivative-free optimizer}
for solving the minimization procedure
\begin{equation}
\label{eq:funct1}
 \min_{z \in \R^d} J(z), \quad J(z) = \frac{1}{2}\|S(z)-u\|^2_{\Gamma} + \frac{\alpha}{2}\|z\|^2,
\end{equation}
where the final term of our objective functional $J$ is a penalty term acting as regularization with regularization parameter $\alpha>0$. This specific penalty term
in \eqref{eq:funct1} corresponds to Tikhonov regularization. Traditionally, to solve \eqref{eq:funct1} one must resort to gradient
methods. EKI instead operates by updating an ensemble of particles $z^{(j)}_n$, where $1 \leq j \leq J$ denotes the ensemble member index, 
using sample covariances which replace the computation of gradients. We define the following sample means and sample covariances 
\begin{align*}
\quad &\bar z_n = \frac1J\sum_{j=1}^J z_n^{(j)}, \quad \bar S_n = \frac1J\sum_{j=1}^J S(z_n^{(j)}),  \\
C_n^{zz} = \frac1J \sum_{j=1}^J (z_n^{(j)}-&\bar z_n)(z_n^{(j)}-\bar z_n)^\top, \quad  C_n^{zp} = \frac1J \sum_{j=1}^J (z_n^{(j)}-\bar z_n)(S(z_n^{(j)})-\bar S_n)^\top\,.  
\end{align*}
The update formulae for EKI are then given as
\begin{align}\label{eq:EKI}
&z_{n+1}^{(j)} = z_n^{(j)} + C_n^{zz} (C_n^{zp} + \Gamma)^{-1}(u_{n+1} ^{(j)}- S(z_n^{(j)})), \\
&u_{n+1}^{(j)} =u + w^{(j)}_{n+1},  \label{eq:obs}
\end{align}
%with defined sample means and sample covariances and
where $(u_{n}^{(j)})$ are perturbed observations with independent realizations of the observational noise $w^{(j)}_{n}\sim\mathcal N(0,\Gamma)$. The above formulation of EKI was originally derived in \cite{ILS13}, and motivated through the ensemble Kalman
filter (EnKF), and its application in reservoir modelling \cite{GE09,GE03,LR09,TMK16}. 
\textcolor{black}{The EnKF is a computational method for the  data assimilation problem. It concerns estimating a stochastic dynamical system $z_n$ with sequential observation $u_n$. In particular, it is of interest to obtain the Bayesian posterior $p(z_n|u_1,\ldots,u_n)$. The inverse problem can be viewed as a special case of data assimilation, where the dynamical system is gorvened by the identity map.
 We refer the reader to the various works on the connections of the EnKF to inverse problems, in the context of EKI \cite{ILS13,LSZ15,LR09,ORL08}.}
\\
Since this seminal work, numerous extensions, and developments, have been made which include deriving regularization schemes, in particular Tikhonov regularization \cite{CST19,MAI16,IY21,WCS22}, deriving analysis both in the continuous and discrete setting \cite{BSW18,BSWW19,BSWW22,CSW19, CT22,SS17,TM2023,W2022} and providing connections with sampling and mean-field
analysis \cite{DL21,GHLS19}. With all these developed works, EKI has thus far only been considered in a framework where the operator $S:\R^d \rightarrow \R^p$ is 
static. As mentioned, in the context of inverse problems, it is an important extension to dynamic forward models where limited literature exists. Our specific interest is the potential connection with the EKI
methodology. dynamic, or time-dependent, inverse problems are a relatively new field of interest, which are motivated through the use of some natural applications.
\textcolor{black}{
 Some of the recent research in this field are through the development of new methodology, such as Kaczmarz-based methods, for linear time-dependent
 problems, as well as developing theory through the use of Bochner spaces, and the connections to well-known models in applied mathematics such as
 the Fokker-Planck equation \cite{AK23,BK17,TTN19}. Time-dependent inverse problems also have connections with ``\textit{all-at-once} inversion
 which is concerned with parameter estimation, related to both the model and observational equation. This notable work has been conducted by Kaltenbacher \cite{BK17,KSW21,KSW21b}. Furthermore, there has been considerable interest in dynamic inverse problems related to imaging applications. These include such tomography-based techniques like 4D computed tomography, positron
emission tomography, or magnetic resonance imaging \cite{HOS21,KMW15,LS17,TR13}. We emphasise that these works consider derivative-based methodologies, which can be cumbersome due to the difficult nature of this dynamic problems. This leaves an open question of whether a computationally sound method, potentially avoiding derivative computations, can be conducted, while also considering different data models. Finally it is worth pointing out that our setting  is different from the standard data assimilation setting. 
We assume the parameter is invariant while the observation operation is dynamical. In most data assimilation problems, it is opposite.  
}
 %\swtd{I still don't see the connection to all-at-once inversion}
 %\swtd{why are time-dependent inverse problems connected to the all at once approach?}
%\swtd{Why are there connections to all-at-once inversion? And if this is true, we should include a reference.} 

\subsection{Problem formulation}
 We consider the task of recovering an underlying \textcolor{black}{ground truth} $z_\ast\in \R^d$, given a time-dependent observation model (in discrete time) of form
\begin{equation}\label{eq:obs_model}
 u_t = S_t z_\ast + w_t,\quad t\in\mathbb{N}\,,
\end{equation}
where $S_t\in\R^{p\times d}$ denotes the %time-dependent, 
dynamic, possibly stochastic\footnote{To be more precise, in case of a stochastic observation operator we understand $S_t\in\R^{p\times d}$ as measurable mappings $S_t:\Omega\to \R^{p\times d}$, i.e.~as random variables taking values in $\R^{p\times d}$.}, %, aperiodic and ergodic 
observation operator mapping to a lower dimensional observation space $\R^p$. {We will assume that the sequence of observation operators $(S_t)_{t\in\mathbb N}$ is either aperiodic and ergodic, or periodic.} For each fixed time $t\in \mathbb{N}$ we assume that the observation operator $S_t$ is linear. Moreover, we assume that the observation is perturbed by independent and identically distributed (i.i.d.) noise $(w_t)_{t\in\N}$ with $w_t\sim\mathcal N(0,\Gamma)$, $t\in\mathbb{N}$, where $\Gamma\in\R^{p\times p}$ is assumed to be symmetric and positive definite. In the following, we denote the natural filtration generated by $(u_t)_{t\in\N}$ as $\mathcal F_t =\sigma(u_s,s\le t)$. Note that $(w_{s})_{s>t}$ is independent of $\mathcal F_t$ by construction. Given a sequence of observations $(u_t)_{t\in\mathbb{N}}$ we want to recover the \textcolor{black}{ground truth} $z_\ast$ sequentially.

In order to discuss the relevance of the considered observation model \eqref{eq:obs_model}, we consider a number of useful examples. Further example are
provided in the following book which is a collection of works on time-dependent inverse problems \cite{KSW21}, %\swtdd{I believe Kaltenbacher is just an editor of this book, which is a collection of papers. Maybe we should cite this in a different style.}, 
such as molecular localization microscopy, acoustic parameter imaging and others.
We now present two model examples, where the first one of them is studied numerically in more details in Section~\ref{sec:num}. %which discusses each particular form of data that we will use, related to our dynamic inverse problem.
%\begin{example}[PDE-constrained optimization]
%\todo[inline]{ToDo: Include examples, in which our assumptions are satisfied.}
%\end{example}

%\todo{Modify this example to one where the unknown is not time-dependent.}
%\swtd{This example is still not included in our observation model. We assume that the unknown $z$ is constant in $t$.}
%\begin{example}[Dynamic Computerized Tomography (DCT)]
%\label{exp:1}
%DCT is concerned \\ with the recovery of some parameter $z$ from X-ray measurements.
%recover the interior $u$ of an object from X-ray measurements \cite{HLR17,HH07}. If
%the object undergoes a motion (e.g. if the patient is moving), then the function $z$ to
%be recovered depends on time and the forward operator is given as
%$$u=S_tR(z_t),$$
%where $R$ denotes a Radon transform and $S_t$ describes the measurement geometry.
%We note that with this particular example both the forward operator and unknown depend
%on time.
%\end{example}
%\swtd{It is good to have examples as below. However, is there any hope to verify our conditions on correlated, iid or periodic data? Neil: Simon I have used iiid now}

\begin{example}[Darcy flow]\label{ex:darcyflow}
We consider the following $2$-dimensional elliptic PDE model
\begin{equation}\label{eq:darcy_flow}
\begin{cases}
-\nabla \cdot (\exp(a)\nabla p) = z^\ast,& x\in D\\
p=0,& x\in\partial D
\end{cases}\, ,
\end{equation}
with domain $D:=(0,1)^2$ and subject to zero Dirichlet boundary conditions. \textcolor{black}{The Darcy flow model is commonly used to model fluid movement in subsurface environments, particularly in groundwater flow. In our case, the coefficient $\exp(a)$ may represent the permeability of the subsurface material, quantifying its ability to transmit fluids. The inverse problem involves recovering the unknown source term $z^\ast\in L^\infty(D)$ from discrete observation points of the pressure $p\in \mathcal V:= H_0^1(D)\cap H^2(D)$ solving equation \eqref{eq:darcy_flow}.}
%Our aim is to recover the unknown source term $z^\ast\in L^\infty(D)$ from discrete observation points of the solution $p\in \mathcal V:= H_0^1(D)\cap H^2(D)$. 
Given $a\in L^\infty(D)$ let $G_a:L^\infty(D)\to \mathcal V$ denote the solution operator of \eqref{eq:darcy_flow} and define an observation operator to be a linear operator $\mathcal O_{[x_1:x_k]}:\mathcal V\to \R^K$ that evaluates $p\in\mathcal V$ in $K$ specified observation points $x_1,\dots,x_K\in D$. More precisely, this means $\mathcal O_{[x_1:x_k]} p = (p(x_1),\dots, p(x_K))^\top\in\R^K$. We are interested in recovering $z^\ast$ given a sequence of observations \[u_t = \mathcal O_t\circ {G_{a_t}}(z^\ast) + w_t\, ,\] where either the solution operator or the observation operator, or both, change in time. Moreover, in each observation step we assume that $w_t\sim\mathcal N(0,\Gamma)$ is independent additive Gaussian noise. In the following, we present a row of specific observation models in which we are interested. 
\begin{itemize}
\item[(i)] Independent observation model: One simple scenario would be the case where $a$ is fixed, but for each measurement the observation operator draws observation points uniformly and independently in $D$. This means, we define $\mathcal O_t := \mathcal O_{[x_1(t):x_K(t)]}$ where $x_i(t)\sim \mathcal U(D)$ are independent for all $i=1,\dots, K$ and $t\ge1$. In this case, we consider a sequence of observations given by
\[u_t = S_t z^\ast + w_t,\]
where $S_t := \mathcal O_t\circ G$ with $G := G_a$ for a fixed $a$.
\item[(ii)] Periodic observation model: Next, we consider the case where the domain $D$ can be decomposed into $p$ disjoint subsets $D_1,\dots,D_p\subset D$ such that $D=\cup_{i=1}^p D_i$. We then define the dynamic observation operator by first selecting a subdomain $j$ and then drawing independently observation points in $D_j$. More precisely, for $t=pk + i$ for $k\ge 0$ and $i=1,\dots p$ we define $\mathcal O_t = \mathcal O_{[x_1(t):x_K(t)]}$ where $x_i(t)\sim\mathcal U(D_i), i\in [K],$ independently. Hence, we move periodically through all sub-domains $D_1,\dots D_p$ and take observations. The sequence of observations is again given by
\[ u_t = S_t z^\ast + w_t,\]
with $S_t := \mathcal O_t\circ G$ and $G:=G_a$ for a fixed $a$.
\item[(iii)] Ergodic observation model: In the last setting, we assume that the observation operator is fixed and given by $\mathcal O := \mathcal O_{[x_1:x_K]}$ for $x_1,\dots,x_K\in\mathcal D$. However, we assume that the diffusion coefficient $a$ is dynamic. Suppose that $(a_t)_{t\ge0}$ is generated from an ergodic Markov chain and define for each state $a_t$ the solution operator ${G_{a_t}}$. This Markov chain may come for example from a Markov chain Monte Carlo algorithm targeting some posterior distribution from a pre-stage Bayesian experiment. In this case, we consider a sequence of observations
\[ u_t = S_t z^\ast + w_t,\]
where this time the dynamic forward model is defined by $S_t := \mathcal O \circ {G_{a_t}}$.
\end{itemize}
We emphasize that this list of examples is by far not complete and the different types of observation models may further be combined. For example, in (ii) one may also consider the case where the sub-domain is picked by some agent who follows a specific decision rule. This may lead again to a limiting ergodic behavior. %Indeed, in our numerical example, we will compare both the periodic and the ergodic choice of the sub-domains. 
\end{example}

In order to reconstruct the underlying ground truth $z_\ast$, the goal is to solve an optimization problem of form
\begin{equation} \label{eq:limiting_Opt}
\min_{z\in\R^d}\ J(z),\quad J(z):= \ell(z-z_\ast) + \frac{\alpha}{2}\|z\|^2,\quad \ell(v):= v^\top A v = \|Sv\|^2\, . 
\end{equation}
{Here $A$ should be interpreted as the log-term average of $S_t^\top S_t$, and $S$ is one of its Cholesky factorization. The exact definition of  $A$ can be found latter in Section \ref{sec:model}.}
$\alpha>0$ denotes a suitable regularization parameter. Note that we assume that $\alpha>0$ has been chosen and is given as fixed parameter. However, we emphasize that the performance of the reconstruction of $z_\ast$ will heavily depend on the choice of $\alpha$. Deriving specific choices of the regularization parameter $\alpha>0$ is left for future work.

The key challenge in the time-dependent setting for solving the optimization problem \eqref{eq:limiting_Opt}, is that the limiting matrix $A$ and the underlying \textcolor{black}{ground truth} $z_\ast$ are unknown, and $z_\ast$ can only be observed through the (random) observation operator $S_t$ perturbed by noise $w_t$, i.e.~through the sequence of observations $(u_t)_{t\in\mathbb N}$ defined in \eqref{eq:obs_model}. Note that the solution of the minimization problem can be written analytically depending on the unknown matrix $A$ and $z_\ast$ as
\[z^{\mathrm{opt}} = (S^\top S+ \alpha I)^{-1}S^\top u_\infty,\quad u_{\infty} := Sz_\ast\,,\]
which solves the first order optimality condition
\[0 = \nabla J(z) = S^\top (Sz-u_\infty) + \alpha z\,. \]

We collect the following notations
\begin{align*}
A = S^\top S, \quad A^{(\alpha)} = S^\top S+ \alpha I,\quad A_t = S_t^\top S_t,\quad A_t^{(\alpha)} =S_t^\top S_t + \alpha I\,,
\end{align*}
and observe that we can write $J(z)=\langle z, A^{(\alpha)} z\rangle$, where $\langle \cdot, \cdot\rangle$ denotes the \textcolor{black}{Euclidean} inner product. Moreover, we have the connection $A^{(\alpha)} z^{\mathrm{opt}} = Az_\ast$, $A_t-A = A_t^{(\alpha)}-A^{(\alpha)}$ and
\[\nabla J(z) = A^{(\alpha)} z - Az_\ast\,. \]
In the upcoming section, we will discuss assumptions on the relation between $S_t$ and $A$ in more details.

\subsection{Model assumptions}
\label{sec:model}
As mentioned, within this work we will consider two forms of data models which are (i) ergodic data, and (ii) periodic data. Note that i.i.d. data can be viewed as special case of both models. To help distinguish each form of data we discuss numerous assumptions related to each data form.  

 \subsubsection*{{Ergodic} data}
 %\swtdd{I would call the first setting "aperiodic" or "ergodic", since both settings give "correlated" data. Neil: Sure i agree.}
 We begin by making the following asymptotic assumption on the dynamic forward operator when using ergodic data.
\begin{assumption}[ergodic data]\label{assum:correlated}
Let $(S_t)_{t\in\mathbb{N}}$ with $S_t\in\R^{p\times d}$ be a sequence of $\mathcal F_t$-adapted (random) matrices such that
%\[\lim_{t\to\infty} \E[S_t^\top S_t]  = S^\top S = A,\]
there exists a symmetric positive definite and bounded matrix $A\in\R^{d\times d}$, i.e.~$0\preceq A \preceq I$, satisfying the following limiting behaviour. For any $\varepsilon>0$ there is an $\tau_\varepsilon \gtrsim |\log(\varepsilon^{-1})|>0$ such that %\sw{($S_t$ assumed to be a Markov process with stationary distribution? What is the correct definition of $A$?)}
%\sw{(For example: For any $\varepsilon>0$ there is an $\tau_\varepsilon>0$ such that
\[ \|\E[S_{t+1+\tau_\varepsilon}^{\top}S_{t+1+\tau_{\varepsilon}}\mid \mathcal F_t]-A\| = \|\E[A_{t+1+\tau_{\varepsilon}}\mid \mathcal F_t]-A\|\le \varepsilon\, ,\]
for all $t\ge1$.
%}
\end{assumption} 
Assumption \ref{assum:correlated} is motivated by the ergodic theory of Markov process \cite{MT93}. In particular, $\tau_\epsilon$ is often interpreted as the mixing time of the process, which describes \textcolor{black}{how fast the process converges} to the stationary distribution. Since the convergence speed is in general exponential, $\tau_\epsilon$ often  scales as $\mathcal{O}(\textcolor{black}{\log \epsilon^{-1}})$. The matrix $A$ in Assumption \ref{assum:correlated}  is the average of $A_t$ under the stationary distribution.

%\swtd{We should discuss this assumption in more details. \textbf{Xin} can you help with this?}

%\swtd{When I am correct, we do not need the assumption $0\prec A$ in all of the settings and indeed positive semi-definite should be enough. This could be used as argue why we do not turn of the regularization in the long-time. Please double check, if I am correct. Xin: if you do regularization, then yes.}
%\todo[inline]{ToDo: discuss this assumption. Is this still some form of Markovian assumption? Or can we say that the assumption is satisfied when $S_t^\top S_t$ is an ergodic Markov chain? Can we say that the assumption is satisfied for periodic observation models?}

\begin{remark}
We emphasize that Assumption~\ref{assum:correlated} is substantially weaker than assuming that the data is generated from an i.i.d. sequence of observations. As result our presented convergence result directly transfer and even simplifies under assumption of i.i.d. data, where we assume that
%\begin{assumption}[i.i.d. data]\label{assum:iid}
%Let 
$(S_t)_{t\in\mathbb{N}}$ %with $S_t\in\R^{p\times d}$ be 
is a sequence of independent and identically distributed (random) matrices such that
\[\E[S_t^\top S_t]  = S^\top S = A,\]
where the expectation matrix $A\in\R^{d\times d}$ is assumed to be positive definite and bounded, i.e.~$0\preceq A \preceq I$.
%\end{assumption} 
\end{remark}
 
 \subsubsection*{Periodic data}
Our second form of data we consider in this work is periodic data. 
%\swtd{Any real-world application example for periodic data in inverse problems?}
%\swtd{skip the following example? I don't think we need it, as we already discussed periodic data in Example 1.1 above}
%\begin{example}
%{
%Let $A \in \mathbb{R}^{2 \times 2}$ be a matrix defined as
%\[ A = \begin{pmatrix}
%a_{11} & a_{12} \\ a_{21} & a_{22}
%\end{pmatrix},
%\]
%which is a (symmetric) strictly positive definite matrix, and consider the periodic system
%\begin{align*} 
%S_{1+pk} &= \begin{pmatrix}
%\sqrt{a_{11}} & 0 \\ 0 & 0
%\end{pmatrix},\quad S_{2+pk} = \begin{pmatrix}
%0 & \sqrt{a_{12}} \\ 0 & 0
%\end{pmatrix},\\ S_{3+pk} &= \begin{pmatrix}
%0 & 0 \\ \sqrt{a_{21}} & 0
%\end{pmatrix},\quad S_{4+pk} = \begin{pmatrix}
%0 & 0 \\ 0 & \sqrt{a_{22}}
%\end{pmatrix},\quad p=4,\ k\ge0\,.
%\end{align*}
%Then for $A_k := S_k S_k^\top$, we have 
%\[\frac1p \sum_{i=n}^p A_{i} = \frac1p A, \]
%\red{[This isn't right. $A_{2+p_k}$ will be diagonal]}
%for all $n\ge1$. Hence, the following condition is always satisfied: For any $\varepsilon>0$, there exists %$\tau_\varepsilon = k\cdot p$ (here even $\varepsilon=0$ and any $k\ge1$), such that
%\[\Big\|\frac1{\tau_{\varepsilon}} A_{n:n+\tau_{\varepsilon}} - \frac1p A\Big\| \le \varepsilon\,. \]
%} 
%\end{example}
Since Assumption \ref{assum:correlated} is known to fail for periodic Markov Chains, we describe %This leads us to 
our next assumption of periodic limiting behavior of the considered time-dynamical observation model. Note that this assumption also allows for randomness in the time-dependent model.
\begin{assumption}[periodic data]
\label{assum:period}
Let $(S_t)_{t\in \mathbb{N}}$ with $S_t\in\R^{p\times d}$ be a sequence of (random) matrices such that
%\[\lim_{t\to\infty} \E[S_t^\top S_t]  = S^\top S = A,\]
%where the 
there exists a positive definite and bounded matrix $A\in\R^{d\times d}$, i.e.~$0\preceq A \preceq I$, satisfying the following limiting behaviour. For any $\varepsilon>0$ there is an $\tau_\varepsilon \gtrsim |\varepsilon^{-1}|>0$ such that %\sw{($S_t$ assumed to be a Markov process with stationary distribution? What is the correct definition of $A$?)}
%\sw{(For example: For any $\varepsilon>0$ there is an $\tau_\varepsilon>0$ such that
\begin{align*}
\Big\|\mathbb E[\frac{1}{\tau_{\varepsilon}}\sum_{k=0}^{\tau_{\varepsilon}-1} S_{t+1+k}^{\top}S_{t+1+k}\mid \mathcal F_t]-A\Big\| &=  \Big\|\frac{1}{\tau_{\varepsilon}}\sum_{k=0}^{\tau_{\varepsilon}-1} \E[A_{t+1+k}\mid \mathcal F_t]-A\Big\|\\&=:  \Big\|\frac{1}{\tau_{\varepsilon}}\E[A_{(t+1):(t+1+\tau_{\varepsilon})}\mid \mathcal F_t]-A\Big\|\le \varepsilon\,,
\end{align*}
for all $t\ge1$.
%}
\end{assumption} 
%\red{For Markov chains with periodicity $p$, $\tau_\epsilon$ can often be chosen as a multiple of $p$. }
%\red{[This assumption is more general than the previous two. So do we need them?]}
%\swtd{What may be $\tau_\varepsilon$ for the periodic assumptions? \textbf{Xin} can you help with this comment?}
Considering periodic data is interesting for a number of reasons, firstly because the extension covers the case where $S_t$ is a periodic sequence,
 whereas the cases of Assumption~\ref{assum:correlated} %and Assumption~\ref{assum:iid}
% \begin{equation}
% \label{tmp:A-A}
%\E_n[A-A_{n+\tau}]\leq \epsilon.
%\end{equation}
 does not cover periodic sequences. %As well as this there are numerous applications for which periodic time-dependent observations can occur. 
%\swtdd{particular examples would be good?} 

 %Particular examples of this
% are medical imaging, such as the Radon transform which is known to have data that resembles a sinusoid. %Another particular example of this is \todo{Neil to finish this!}
%\swtd{if this is true, can we use this as numerical example? Or can we at least state the model as example. Right now, above we have two examples, where we can not see whether one of our assumptions is satisfied.}

%\todo{Extend this to a weaker assumption with expectation} 
%\swtd{right now this assumption is deterministic wrt $A_k$, we may incorporate an (conditional) expectation again. However, for periodic examples I am not sure if it is needed. On the other side, we may also include iid data under this assumption as well. In principle, it would say, that the MC estimation of $A$ through the average over $A_k$ needs to be sufficiently small. We may use
%\[\E_n\Big[\Big\|\frac{1}{\tau_{\varepsilon}}A_{n:n+\tau_{\varepsilon}}-A\Big\|^2\Big]\le \varepsilon, \]
%which comes from MC for iid data $A_k\sim A$. Connect this also to subsampling EKI \cite{HLS23}.}

%\begin{figure}[h!]
%\centering
%\includegraphics[scale=0.75]{data.pdf}
%\caption{Different representations of the data. The data is represented by blue, and pointwise observations are highlighted in red.}
%\label{fig:data}
%\end{figure} 

\subsection{Our Contributions:}
To conclude this section, we summarize our contributions below.
\begin{itemize}
\item We present a formulation of the ensemble Kalman inversion based on the time-dependent observation model \eqref{eq:obs_model}.
This differs from the conventional EKI, where the forward operator $S$ is usually assumed to be fixed and is applied as tool in a static inverse problems setting. Our resulting
algorithm is entitled ``dynamic EKI". %\swtdd{Markov chain EKI does not fit well? Maybe something like "dynamic EKI"?}
%which is based on Markov chain gradient descent. 
\item We provide a number of theoretical results to demonstrate the validity of our proposed scheme. Our initial analysis
require the controlling of the modified covariance matrix, for which we demonstrate this through lower and upper bounds on the ensemble collapse,
and providing moment bounds. %bounding the first and second moments. 
\item A convergence analysis is provided for the  \textcolor{black}{DEKI} as stochastic optimization method of \eqref{eq:limiting_Opt}. %, which requires a %number of standard assumptions from
%optimization theory, and related to the choice of the 
%a sufficiently small step-size $h>0$. Our analysis  
This analysis covers the three different data types of interest: (i) i.i.d. data, (ii) ergodic data and (iii) periodic data.
\item Numerical experiments are conducted verifying the theory that is attained. We test our \textcolor{black}{DEKI} on a 2D Darcy flow PDE example, comparing our ergodic, periodic and i.i.d. data.
\end{itemize}

\subsection{Outline}
The outline of this work is as follows: In Section \ref{sec:EKI} we describe and provide our dynamic version of
EKI. %We also discuss a number of necessary definitions, which we will make use of and state our main convergence
%theorem. 
This will lead onto Section \ref{sec:bounds} which is where our preliminary analysis is presented. %For this
%analysis we require the controlling of our newly defined sample covariances, for which we are able to show an ensemble
%collapse result, but also moment bounds. 
Section \ref{sec:main} is devoted to the proof of our theorems, which are seperated based on the type of data that is used. Numerical results are shown in Section \ref{sec:num} %on a PDE example comparing our data-types, verifying our theory. 
and finally, we conclude our findings in Section \ref{sec:conc}.

\section{Ensemble Kalman inversion with dynamic forward operator}\label{sec:EKI}
In this section we introduce and discuss our proposed algorithm, referred to as \textcolor{black}{dynamic-EKI (DEKI)}.
In order to derive it, we firstly present the vanilla version of EKI \eqref{eq:EKI} in a modified setting, which 
provides strong differences. % We also introduce some definitions that are required later, and we finish
We finish this section by stating our main result which is a convergence theorem for \textcolor{black}{DEKI} with
given rate of convergence. 

%We begin by introduce some assumptions which we utilize for our theory.
%\begin{defn}[$L$-smooth]
%\label{def:Ls}
%The functional $J:\R^d\to\R$ is said to be $L$-smooth if it is differentiable with $L$-Lipschitz continuous gradient, i.e.~there exists $L>0$ such that
%$$
%\| \nabla J( z) -\nabla  J( z')\| \leq L\| z-z'\|, \quad z,z'\in\R^d\, .
%$$
%\end{defn}
%\begin{defn}[strong-convexity]
%\label{def:sc}
%A continuously differentiable functional $J:\R^d\to\R$ is said to be $\alpha$-strongly convex if the following bound holds
%$$
%\langle \nabla J( z) - \nabla J( z' ) ,  z- z'\rangle  \geq \alpha \| z- z'\|^2,\quad z,z'\in\R^d ,
%$$
%where $\alpha>0$.
%\end{defn}
%As a direct consequence of $\alpha$-strong convexity and $L$-smoothness, we have the Polyak-Lojasiewicz (PL) inequality, or condition given below.
%\begin{defn}[Polyak-Lojasiewicz condition]
%\label{def:PL}
%Let us assume that the differentiable functional $J:\R^d\to\R$, \eqref{eq:funct}, is $\alpha$-strongly convex and that
%the unique minmizer is denoted as $z^{\mathrm{opt}}$
%inequality hold
%$$
%J(z^{\mathrm{opt}}) \geq J(z) - \frac{1}{2\alpha}\| \nabla J((z)\|^2,\quad z\in\R^d.
%$$
%\end{defn}

One possible way to solve \eqref{eq:limiting_Opt}, is to apply a gradient descent scheme. Since we assume that there is no access to $A$, $S$ and $z_\ast$, we may formulate the dynamic (stochastic) gradient descent scheme by
\begin{equation}\label{eq:gradientdescent}
z_{t+1} = z_t - \eta_t \nabla J_{t+1}(z_t) = z_t - \eta_t \left(S_{t+1}^\top (S_{t+1} z_t - u_{t+1}) + \alpha z_t\right)\,, 
\end{equation}
with initial $z_0\in\R^d$, $t\in\mathbb N$, where $(\eta_t)_{t\in\mathbb{N}}$ with $\eta_t>0$, $t\in\mathbb{N}$ denotes a sequence of step sizes, and $J_t:\R^d\to\R$ denotes the time-dependent loss function defined by
\begin{equation}
\label{eq:funct}
J_t(v) = \|S_t v - u_t\|^2 + \frac{\alpha}2\| v\|^2,\quad t\in \mathbb{N}\,.
\end{equation}

%\todo[inline]{ToDo: Motivate why this can be seen as stochastic optimization scheme. Connect to Markov Chain gradient descent?}

We propose to apply an alternative algorithm motivated by the EKI, which can be written in simplified form as %\swtd[inline]{This corresponds to a simple preconditioned gradient descent scheme. Do we want to ignore the formulation including the Kalman gain?}
\begin{equation}\label{eq:unpert_EKI}
z_{t+1}^{(j)} = z_t^{(j)} - \eta_t C_t^{zz} \left(S_{t+1}^\top (S_{t+1} z_t^{(j)} - u_{t+1}) + \alpha z_t^{(j)}\right),
\end{equation}
where our sample covariance and mean are defined as
$$
C_t^{zz} = \frac1J \sum_{j=1}^J (z_t^{(j)}-\bar z_t)(z_t^{(j)}-\bar z_t)^\top, \quad \bar z_t = \frac1J\sum_{j=1}^J z_t^{(j)}.
$$ 
Note that an alternative formulation of the ensemble Kalman inversion incorporates perturbed observations such that the iteration can be written as 
\begin{equation}\label{eq:pert_EKI}
\begin{split}
z_{t+1}^{(j)} &= z_t^{(j)} - \eta_t C_t^{zz} \left(S_{t+1}^\top (S_{t+1} z_t^{(j)} - u_{t+1}^{(j)}) + \alpha z_t^{(j)}\right),\\
u_{t+1}^{(j)} &= S_{t+1} z_\ast + w_{t+1}^{(j)},
\end{split}
\end{equation}
where $w_t^{(j)}\sim\mathcal N(0,\Gamma)$ are assumed to be independent realizations of the noise. 
Based on this, we present our new methodology in Algorithm \ref{alg:EKI_cor}. Note that the interacting particle system generated by \eqref{eq:unpert_EKI} or \eqref{eq:pert_EKI} respectively are adapted with respect to the filtration $(\mathcal F_t)_{t\in\N}$ by construction. In our convergence results %presented in Section~\ref{ssec:main_results} 
we focus on the \textcolor{black}{DEKI} with unperturbed observations \eqref{eq:unpert_EKI}.

\begin{algorithm}[htb!]
\begin{algorithmic}[1]
\State \textbf{Input:} \begin{itemize}
 \item initial ensemble $z_0^{(i)}\in \R^d$, $i=1,\dots,J$,
 \item sequence step sizes $(\eta_t)_{t\in\N}$.
 \end{itemize}
\For{$t=1,\dots, T$}
	\State perturb observations $u_{t+1}^{(j)} \sim \mathcal N(u_{t+1},\eta_{t}^{-1}\Gamma)$,
	\State for $j=1,\dots,J$ iterate 
	\[
	z_{t+1}^{(j)} = z_t^{(j)} - \eta_t C_t^{zz} \left(S_{t+1}^\top (S_{t+1} z_t^{(j)} - u_{t+1}^{(j)}) + \alpha z_t^{(j)}\right)\, ,
	\]
	\State set $\bar z_t = \frac1J \sum_{i=1}^J z_t^{(i)}$.
\EndFor
\end{algorithmic}
 \caption{Dynamic Ensemble Kalman Inversion (DEKI)}
   \label{alg:EKI_cor}
\end{algorithm}

%\swtd{Oh no, I did a fundamental mistake: We have to shift the time index in $S_t$ in \eqref{eq:unpert_EKI} to 
%\[z_{t+1}^{(j)} = z_t^{(j)} - \eta_t C_t^{zz} \left(S_{t+1}^\top (S_{t+1} z_t^{(j)} - u_{t+1}) + \alpha z_t^{(j)}\right),\]
%such that $z_t$ is $\mathcal F_t$-measurable. I have to do this throughout the proofs below. This does not change anything from the theoretical analysis, it is just a technical detail which is needed when taking the conditional expectations $\E[\cdot \mid \mathcal F_t]$. I will do this after lunch time.
%}

The above formulation for our \textcolor{black}{DEKI} may seem quite different to its original form,
presented in \eqref{eq:EKI}-\eqref{eq:obs}. To help to understand the intuition behind it we briefly
discuss and present the connection to gradient related algorithms which follow similar ideas.
%\swtd{Discuss SGD, MCGD, but also (block-) coordinate descent}
%\paragraph{\it Stochastic gradient descent}
Assuming that $S_t^\top S_t$ is an independent and unbiased estimator of $S^\top S = A$ for each $t\in\mathbb N$, we can view \eqref{eq:gradientdescent} as specific form of stochastic gradient descent (SGD), where for each $t\in\N$ one can verify 
\[\mathbb E[\nabla J_{t+1}(z_t) \mid\mathcal F_t] = \nabla J(z_t)\,.\]
As result, our proposed \textcolor{black}{DEKI} can be viewed as preconditioned SGD algorithm. Indeed, under this assumption, we may view \eqref{eq:unpert_EKI} as discrete time variant of 
%\\\\
%By relaxing Assumption \ref{assum:correlated}, this indicates we are now in the setting of i.i.d. data.
%The correlated data example, related to EKI, goes beyond the use of i.i.d. data which has been connected with EKI
%in a recent paper 
of the subsampling approach for EKI recently proposed by Hanu et al. \cite{HLS23}. In that work the authors provided a subsampling approach to EKI  in the continuous-time formulation, where the different observations are chosen using switching times.
%from the use randomization methods such as stochastic gradient descent in optimization. 
It is worthy to point out that no particular
 rate of convergence was derived, as a continuous-time setting was primarily adopted. As result we can view the i.i.d. setting as special case of our assumption of ergodic data Assumption~\ref{assum:correlated} for which we derive a rate of convergence in Corollary~\ref{thm:main3}. %Therefore this raises the question of 
 %whether we can attain a rate, which could depend on the correlated data setting.
%\paragraph{\it Markov chain gradient descent}
%\swtd{change to updated notation}
%, a similar algorithm, known as Markov chain gradient descent (MCGD). 
Under Assumption~\ref{assum:correlated} the iterative scheme \eqref{eq:gradientdescent} is related to Markov chain gradient descent (MCGD). This method operates in a similar fashion to SGD, with the exception of using data coming from an ergodic Markov chain. 
%Let us assume we have correlated, then the update formula for MCGD
%is
%$$
%z_{t+1} = z_t- \eta_{t+1}\big((S^{\bsy_t})^\top(S^{\bsy_t}z_t - u_0) + \alpha z_t \big).
%$$
%$$
%z_{t+1} = z_t- \eta_{t+1}\big(S_{t}^\top(S_{t}z_t - u_0) + \alpha z_t \big).
%$$
This algorithm has been of particular interest within the machine learning community \cite{ME23,SL19,WLY23}. 
%\begin{algorithm}[htb!]
%\begin{algorithmic}[1]
%\State \textbf{Input:} \begin{itemize}
% \item sequence of step sizes $(\eta_t,t\in\N)$,
% \item initial data distribution $\mu_0$,
% \item $(u_t,t\in\N)$ generated by \eqref{eq:obs_model}
% \item initial state $z_0\in L^2(D)$.
% \end{itemize}
% \State \textbf{Output:} 
%\For{$t=1,\dots, T$}
%	%\State compute $g(z_t,u_t) =$,
%	\State iterate $z_{t+1} = z_t - \eta_{t} \left(S_{t+1}^\top(S_{t+1}z_t - u_{t+1}) + \alpha z_t\right)$.
%\EndFor
%\end{algorithmic}
% \caption{Markov Chain Gradient Descent Method}
% \label{alg:MCGD}
%\end{algorithm}
%\swtd[inline]{I would be a bit careful to call this scheme "Markov chain gradient descent", as this is not the particular form in which MCGD is formulated. The ideas are related, and proofs from MCGD may transfer to this setting, but it is a different algorithm just from the update formula.\\
%Suggestion: We don't make the connection too prominent through writing it as own subsection. We just write it as remark and do not use an own algorithm environment. }
%\paragraph{\bf Coordinate descent method}
%Under Assumption~\ref{assum:period}\dots

%\subsection{Properties of EKI for linear operators}
Before we continue with stating %various assumptions, and 
our main theorems %, we require on our functional 
we briefly 
state the motivation behind the use of EKI-based methodology. 
A natural question to ask is why to consider EKI, compared to other well-known methods such as gradient, or stochastic gradient descent.
There are three main reasons which constitute to our motivation in EKI. From a technical discussion,
the preconditioning that EKI attains through the covariance $C^{zz}_t$ results in several advantages. One of them (i) is that it will satisfy the affine
invariance property. This would be similar to methods which require Hessian information such as Newton-type methods. Furthermore, (ii) we can potentially 
make use of the preconditioner. In particular, it can be treated as a way to control the learning rate, where below we set $\eta_t=h>0$ sufficiently small but fixed. Finally, (iii) the implementation of 
Algorithm \ref{alg:EKI_cor} avoids the computation of $S_t^\top$ by computing the cross covariance $C^{zS}$ instead of $C_t^{zz}S_t^\top$. This alternative computation
can save on associated computational cost, especially for higher dimensional problems. 

In the remaining manuscript, we will focus on EKI with unperturbed observation, i.e.~on the convergence behavior of $(\{z_t^{(j)}\}_{j=1}^J, t\in\mathbb{N})$ generated by \eqref{eq:unpert_EKI}. 

\subsection{Main result}\label{ssec:main_results}
In this section, we state our main result of this article which is a convergence result with given rate. The convergence is quantified through the expected loss %of the difference of functionals evaluated between the optimization
%solution and  
evaluated at the ensemble mean. 
In order for us to derive the main theorems we require a controllability of our 
covariance $C^{zz}_t$, which can be controlled through the ensemble spread
defined 
\begin{equation}
\label{eq:spread}
e_{t}^{(j)} = z_{t}^{(j)} - \bar z_{t},\quad t\ge0.
\end{equation}
Our first main result presents the convergence under Assumption~\ref{assum:correlated}.
\begin{thm}[ergodic data]\label{thm:main}
Let $(\{z_t^{(j)}\}_{j=1}^J, t\in\mathbb{N})$ be generated by \eqref{eq:unpert_EKI} with fixed $\eta_t=h>0$ and %(almost surely) linearly independent 
initial ensemble $\{z_0^{(j)}\}_{j=1}^J$, such that $C_0^{zz}\succ \sigma_l I$ almost surely for some $\sigma_l>0$ with $\lambda= \sigma_l q\in(0,1)$ and $q=2\mu$. Moreover, let 
$$
h\le \min\Big(\frac{\alpha}{J(A_{\max}+\alpha)^2},\frac{J}{\alpha}\Big) E_0^{-1},
$$ 
where $E_0$ is defined as $E_0=\frac{1}{J}\sum_{j=1}^J \|e_0^{(j)}\|^2$, from \eqref{eq:spread}.
Under Assumption~\ref{assum:correlated} for any $\varepsilon>0$ we have that 
%\[\E[ J(\bar z_t) - J(z^{\mathrm{opt}})] \in \mathcal O\left(\left(\frac{1}{T}\right)^{\lambda} (1+\tau(\varepsilon)) + \sum_{t=0}^T \left(\frac{t+1}{T}\right)^\lambda \frac{1}{(t+1)^2} + \varepsilon\right)\,, \]
\begin{align*}
%\E[ J(\bar z_t) - J(z^{\mathrm{opt}})] \in \mathcal O\Bigg(& \sw{\varepsilon+\left(\frac{1}{T}\right)^{\lambda} ((1+\log(\tau_\varepsilon))^2+\log(T))}\\&\quad \sw{+ (1+\log(\tau_\varepsilon)) (1+\log(T)) \sum_{t=0}^T \left(\frac{t+1}{T}\right)^\lambda \frac{1}{(t+1)^2} }\Bigg),
\E[ J(\bar z_t) - J(z^{\mathrm{opt}})] \in \mathcal O\Bigg(&\varepsilon+\frac{1+\tau_\varepsilon (1+\log(\tau_\varepsilon)+\log(T))}{T^{\lambda}} %\\ &\sw{+ {\tau_\varepsilon (1+\log(T))} \sum_{t=0}^T \left(\frac{t+1}{T}\right)^\lambda \frac{1}{(t+1)^2} }
\Bigg) \,.
\end{align*}
%\red{$\epsilon$ should be placed earlier, like first line? Else it is part of the summation? Same goes for all other results}
%where $\lambda\ge \frac{1}2 \sigma_l q$ and $q=2\mu$.
In particular, if we pick $\varepsilon = \frac{1}{T}$, resulting in $\tau_\varepsilon=\log(T)$, we obtain
%\begin{align*}
%\E[ J(\bar z_t) - J(z^{\mathrm{opt}})] \in \mathcal O\Bigg(&\left(\frac{1}{T}\right)^{\lambda} {(1+\log(T)+\log(\log(T))+\log(T)^2)}\\ &+ {(1+\log(T)+\log(T)^2)}\sum_{t=0}^T \left(\frac{t+1}{T}\right)^\lambda \frac{1}{(t+1)^2} + \frac1T \Bigg)\,.
%\end{align*}
\begin{equation*}
\E[ J(\bar z_t) - J(z^{\mathrm{opt}})] \in \mathcal O\left( \frac{\log(T)^2}{T^{\lambda}} \right)\,.
\end{equation*}
%\red{The result can be simplified? like $\log T^2>>\log T>>log log T$. The summation can be bounded by the other two terms, so should be removed? Should do similar things for next result?}
\end{thm}

Let us mention again, that the convergence analysis simplifies if we replace the ergodic data assumption by assuming an i.i.d.~observation model.
As mentioned in Section \ref{sec:intro} the work of Hanu et al. \cite{HLS23}, provide a subsampling approach to EKI which is related to our \textcolor{black}{DEKI} in discrete time using i.i.d. data. %However their work does not provide a rate of convergence. %hence
%why 
We provide the following corollary which attains a convergence rate in this setting.
%\swtd{Discuss the subsampling approach in Hanu et al in more detail in the introduction and emphasize that we obtain a convergence rate for the discrete EKI here. When I am correct, they don't have a rate of convergence}
\begin{cor}[i.i.d. data]\label{thm:main3}
Let $(\{z_t^{(j)}\}_{j=1}^J, t\in\mathbb{N})$ be generated by \eqref{eq:unpert_EKI} with fixed $\eta_t=h>0$ and %(almost surely) linearly independent 
initial ensemble $\{z_0^{(j)}\}_{j=1}^J$, such that $C_0^{zz}\succ \sigma_l I$ almost surely for some $\sigma_l>0$ with $\lambda= \sigma_l q\in(0,1]$ and $q=2\mu$. Moreover, let 
$$
h\le \min\Big(\frac{\alpha}{J(A_{\max}+\alpha)^2},\frac{J}{\alpha}\Big) E_0^{-1},
$$ 
where $E_0$ is defined as $E_0=\frac{1}{J}\sum_{j=1}^J \|e_0^{(j)}\|^2$, from \eqref{eq:spread}.
Assume that the sequence of matrices $(S_t)_{t\in\mathbb N}$ are i.i.d. with expectation $\E[S_1^\top S_1]= A$. Then we have that 
\begin{align*}
\E[ J(\bar z_t) - J(z^{\mathrm{opt}})] \in \mathcal O\left(\frac{1}{T^{\lambda}} %+\sum_{t=0}^T \left(\frac{t+1}{T}\right)^\lambda \frac{1}{(t+1)^2}
\right)\,.
\end{align*}
%where $\lambda\ge \frac{1}2 \sigma_l q$ and $q=2\mu$.
\end{cor}
Comparing the convergence results presented in Theorem~\ref{thm:main} and Corollary~\ref{thm:main3}, we observe that we achieve nearly the same asymptotic convergence behavior under ergodic data as compared to i. data (up to an additional factor of $\log(T)^2$).

Next, we consider the extension from the above result based on periodic data,
which was discussed in Section \ref{sec:intro}. This result is given below, which 
follows very similarly to that of Theorem \ref{thm:main}.
%Preferably, we want \swtd{ToDo: update notation to the one used in Assumption \ref{assum:period}. And in general, move this discussion to Section 1.3}
%\todo{Neil: will do this!}
%\begin{equation}
%\label{tmp:A-A2}
%\E_n[A-\frac{1}{T}A_{n:n+T}]\leq \epsilon,\quad  A_{n:n+T}=\sum_{k=n+1}^{n+T} A_k.
%\end{equation}
%This is an interesting extension simply because it covers the case where $A_n$ is a periodic sequence,
% whereas the case of
 %\begin{equation}
% \label{tmp:A-A}
%\E_n[A-A_{n+\tau}]\leq \epsilon.
%\end{equation}
% does not cover that case. This result is given below. 

%\swtd{include final error bound from the proof}
\begin{thm}[periodic data]\label{thm:main2}
Let $(\{z_t^{(j)}\}_{j=1}^J, t\in\mathbb{N})$ be generated by \eqref{eq:unpert_EKI} with fixed $\eta_t=h>0$ and %(almost surely) linearly independent 
initial ensemble $\{z_0^{(j)}\}_{j=1}^J$, such that $C_0^{zz}\succ \sigma_l I$ almost surely for some $\sigma_l>0$ with $\lambda= \sigma_l q\in(0,1)$ and $q=2\mu$. Moreover, let 
$$
h\le \min\Big(\frac{\alpha}{J(A_{\max}+\alpha)^2},\frac{J}{\alpha}\Big) E_0^{-1},
$$ 
where $E_0$ is defined as $E_0=\frac{1}{J}\sum_{j=1}^J \|e_0^{(j)}\|^2$, from \eqref{eq:spread}.
Then under Assumption \ref{assum:period}, for any $\varepsilon>0$ we have that 
\begin{align*}
%\E[ J(\bar z_t) - J(z^{\mathrm{opt}})] \in \mathcal O\Bigg(&\left(\frac{1}{T}\right)^{\lambda}{(1+(\tau_\varepsilon (1+\log(\tau_\varepsilon)+\log(T)))} \\ &+ {(1+\tau_\varepsilon (1+\log(T)))} \sum_{t=0}^T \left(\frac{t+1}{T}\right)^\lambda \frac{1}{(t+1)^2} + \varepsilon\Bigg)\,,
\E[ J(\bar z_t) - J(z^{\mathrm{opt}})] \in \mathcal O\Bigg(\varepsilon+\frac{1+(1+\log(\tau_\varepsilon))^2+(1+\log(\tau_\varepsilon))(1+\log(T))}{T^\lambda}%\\&\quad \sw{+ (1+\log(\tau_\varepsilon)) \log(T) \sum_{t=0}^T \left(\frac{t+1}{T}\right)^\lambda \frac{1}{(t+1)^2} }
\Bigg).
\end{align*}
%\swtdd{to be changed to $\tau(\varepsilon) = \varepsilon^{-1}$. @Xin: what would be a suitable choice for $\varepsilon$? Maybe $\varepsilon = \frac{1}{T^{\lambda/2}}$, such that we would get $\frac{1}{T^{\lambda/2}}$ with additional $\log$-terms.}
%where $\lambda\ge \frac{1}2 \sigma_l q$ and $q=2\alpha$.
In particular, if we pick $\varepsilon = \frac1T$, resulting in $\tau(\varepsilon)=T$, we obtain
\begin{equation*}
\E[ J(\bar z_t) - J(z^{\mathrm{opt}})] %\in \mathcal O\Bigg(&\left(\frac{1}{T}\right)^{\lambda}{(1+(\mathsf{C} (1+\log(\mathsf{C})+\log(T)))} \\ &+ {(1+\mathsf{C} (1+\log(T)))} \sum_{t=0}^T \left(\frac{t+1}{T}\right)^\lambda \frac{1}{(t+1)^2}\Bigg)\,,
\in \mathcal O\Bigg(\frac{\log(T)^2}{T^{\lambda}} %+ \sum_{t=0}^T \left(\frac{t+1}{T}\right)^\lambda \frac{1}{(t+1)^2}
\Bigg).
\end{equation*}
\end{thm}
Similarly to the ergodic data setting, we achieve nearly the same asymptotic convergence behavior as in the i.i.d. data case.
%\swtd{I would skip this table}
%We provide a table to indicate the difference choices of data, and the choices of $\tau_\varepsilon$
%and $\varepsilon$ related to Theorem \ref{thm:main}. This is given below in Table \ref{table:comp}.
%\begin{table}[h!]
%\begin{center}
%\begin{tabular}{ |c|c|c|c| } 
%\hline
%\textbf{Type of Data} & $\varepsilon$ & $\tau_\varepsilon$ \\
%\hline
%{Correlated data $u_t$} & 1/$T$ & $\log T$ \\ 
%Periodic data $u_t$ & 0 & $\mathsf{C}$ \\ 
%i.i.d. data $u$ & ? & 0  \\ 
%\hline
%\end{tabular}
%\end{center}
%\caption{Comparison of our different forms of data in terms of $\varepsilon$ and $\tau_\varepsilon$. For the periodic case, we denote $\mathsf{C}>0$ to be some constant.}
%\label{table:comp}
%\end{table}
We will defer the proof, of all the above states results to Section \ref{sec:main}.
%\red{Please do the todo note.} 
%\swtd{There was a mistake in translating the final error bound from the proof to the Theorem. We can again simply take $\varepsilon=1/T$.}

\textcolor{black}{
\begin{remark}
\label{rem:ensemble}
It is worth noting that for $C_0^{zz}\succ \sigma_l I$ to hold, it is necessary that the ensemble size satisfies $J>d+1$. The necessity can partly be explained by the subspace property, explained later in Lemma \ref{lem:subspace}, as otherwise DEKI will stay in a subspace defined by the initial ensemble. This can be a computational issue for high dimensional problem. Localization in \cite{TM2023} and dropout techniques \cite{LRT23} may be implemented on DEKI to resolve this issue. 
% An important question, related to the above theorems, is how sensitive they are with respect to the choice of the ensemble size $J$
% and the model dimension $d$? This is an important question related to the ensemble Kalman filter in general. Our Theorems do not discuss such considerations of $J<d$. For our numerical experiments we will consider a choice of $J$ not very small but not large, but with the constraint of $J>d$. To consider the other case one could exploit more sophisticated methods such as localization. This would require a careful construction, which could depend on the data-type. We believe this is extendable, but goes beyond the scope of this work and leave this as future work.
\end{remark}}

\section{Ensemble collapse and moment bounds}
\label{sec:bounds}

In this section, we present our preliminary analysis required before proving Theorem \ref{thm:main}.
Our two results which we utilize include an ensemble collapse result, stating that our ensemble of particles will collapse to a single point but not too fast. This is required to bound our sample covariance matrix from above and below. Our second result in this section
is a moment bound related to the evolution of the dynamical system \eqref{eq:unpert_EKI}.
\\
\textcolor{black}{We begin by discussing the invariant subspace property of DEKI. This is inherited from EKI \cite{ILS13} and its  regularized versions: 
% The  states that the solution of EKI is spanned by the linear span of the initial ensemble. In the context of DEKI this also holds, independent on the choice of data discussed above. We provide a lemma below which states this. 
\begin{lem}[subspace property]
\label{lem:subspace}
If $\mathcal{S}$ is the linear subspace spanned by the initial ensemble $\{z^{(i)}_0\}_{i=1}^{J}$, then $z^{(j)}_t \in \mathcal{S}$ for all $(t,j) \in \mathbb{Z}^{+} \times \{1,\ldots,J\}$.
\end{lem}
\begin{proof}
The proof is identical to that presented in \cite{CST19} and \cite{ILS13}. In particular, recall the update \eqref{eq:unpert_EKI}:
\[
z_{t+1}^{(j)} = z_t^{(j)} - \eta_t C_t^{zz} \left(S_{t+1}^\top (S_{t+1} z_t^{(j)} - u_{t+1}) + \alpha z_t^{(j)}\right).
\]
Using mathematical induction, it is easy to check that  if $z_t^{(j)}\in \mathcal{S}$, the range of operator $C^{zz}_t$ will also be inside $\mathcal{S}$, and $z_t^{(j+1)}\in \mathcal{S}$. 
\end{proof}
}
We now start the discussion with the control of the empirical covariance matrix $C_t^{zz}$. The iterative evolution of the ensemble mean $\bar z_t$ is given by
\[\bar z_{t+1} = \bar z_t -  hC_t^{zz}\left( S_{t+1}^\top (S_{t+1} \bar z_t - u_{t+1}) + \alpha \bar z_t\right),  \]
and the ensemble spread $e_t^{(j)} = z_t^{(j)} - \bar z_t$ follows the iteration
\begin{align*}
e_{t+1}^{(j)} = z_{t+1}^{(j)} - \bar z_{t+1} &=  z_{t}^{(j)} - \bar z_{t} -  C_t^{zz} \left(S_{t+1}^\top (S_{t+1} (z_t^{(j)}-\bar z_t))+ \alpha (z_t^{(j)}-\bar z_t)\right)\\
& = e_t^{(j)} - h C_t^{zz} \left(S_{t+1}^\top S_{t+1} e_t^{(j)} - \alpha e_t^{(j)}\right).
\end{align*}
By definition of the sample covariance matrix $C_{t+1}^{zz}=\frac1J\sum_{j=1}^J e_{t+1}^{(j)}( e_{t+1}^{(j)})^\top $, we have 
\begin{align*}
C_{t+1}^{zz} %&= \frac1J\sum_{j=1}^J e_{t+1}^{(j)}( e_{t+1}^{(j)})^\top\\ 
&= C_t^{zz} - 2C_t^{zz} (S_{t+1} S_{t+1}^\top + \alpha I)C_t^{zz} \\ &\quad +C_t^{zz} (S_{t+1} S_{t+1}^\top +\alpha I) C_t^{zz} (S_{t+1} S_{t+1}^\top +\alpha I) C_t^{zz}\\
&= C_t^{zz} - 2h C_t^{zz} A_{t+1} C_t^{zz} + h^2 C_t^{zz}A_{t+1} C_t^{zz} A_{t+1} C_t^{zz},
\end{align*}
where we have defined $A_t := (S_t S_t^\top +\alpha I)$. 
For the upper bound on the sample covariance, we consider the dynamical evolution of $E_t:=\frac{1}{J}\sum_{j=1}^J \|e_t^{(j)}\|^2$. Note that
\[\|C_t^{zz}\|_{\mathcal F}^2 = \frac{1}{J^2}\sum_{j,l=1}^J \langle e_t^{(j)},e_t^{(l)}\rangle ^2 \ge \frac1J E_t^2, \]
where we have dropped positive terms in the sum over $j,l$ and applied Jensen's inequality. Moreover, we have
\[\|C_t^{zz}\|_{\mathcal F}^2 \le \frac{1}{J^2}\sum_{j,l=1}^J \|e_t^{(j)}\|^2\|e_t^{(l)}\|^2 = E_t^2,\]
by Cauchy-Schwarz inequality.
\begin{lem}[Ensemble collapse]\label{lem:ensemble_collapse}
Let $(\{z_t^{(j)}\}_{j=1}^J, t\in\mathbb{N})$ be generated by \eqref{eq:unpert_EKI} with fixed $\eta_t=h>0$ and %(almost surely) linearly independent 
initial ensemble $\{z_0^{(j)}\}_{j=1}^J$, such that $C_0^{zz}\succ \sigma_l I$ almost surely for some $\sigma_l>0$. Moreover, let 
$$
h\le \min\Bigg(\frac{\alpha}{J(A_{\max}+\alpha)^2},\frac{J}{\alpha}\Bigg) E_0^{-1}.
$$ 
Then for all $t\ge0$ it holds true that \[ \frac{\sigma_l}{h(t+1)} I \preceq C_t^{zz} \preceq \frac{\sigma_u}{(t+1)} I = \frac{J}{h\alpha(t+1)} I\,,\]
almost surely for $\sigma_u = \frac{J}{h\alpha}$.
Moreover, we have that
\[\|C_t^{zz}-C_{t+1}^{zz}\|\le \Big((A_{\max}+\alpha)+h(A_{\max}+\alpha)^2\Big)\Bigg(\frac{J}{\alpha(t+1)^2}\Bigg)\,, \] %\frac{3}{h(t+1)^2}\,. \]
almost surely.
\end{lem}
\begin{proof}
For fixed $j\in\{1,\dots,J\}$ the evolution of $E_t$ is given by
\begin{align*}
E_{t+1} = \frac1J \sum_{j=1}^J\|e_{t+1}^{(j)}\|^2 &= E_t - 2h \frac{1}{J}\sum_{j=1}^J \langle e_t^{(j)} , C_t^{zz} A_{t+1}^{(\alpha)}e_t^{(j)}\rangle  + h^2 \|C_t^{zz} A_{t+1}^{(\alpha)}e_t^{(j)}\|^2\\
 &\le E_t - 2h\alpha \|C_t\|_{\mathcal F}^2 + h^2 (A_{\max}+\alpha)^2 \|C_t\|_{\mathcal F}^3\\
 &\le E_t - \frac{2h\alpha}{J} E_t^2 + h^2 (A_{\max}+\alpha)^2 E_t^3.
\end{align*}
Suppose that $h\le \Big(\frac{\alpha}{J(A_{\max}+\alpha)^2}\Big) E_0^{-1}$, then it follows that
$h^2 (A_{\max}+\alpha)^2 E_0^3 \le \frac{h\alpha}{J} E_0^2$, implying that $E_1\le E_0$. This again implies $h\le \big(\frac{\alpha}{J(A_{\max}+\alpha)^2}\big) E_1^{-1}$ such that $E_2\le E_1$. By induction we hence obtain that $E_t$ is monotonically decreasing and 
\begin{equation}\label{eq:decrease}
h^2 (A_{\max}+\alpha)^2 E_t^3 \le \frac{h\alpha}{J} E_t^2 ,
\end{equation} 
for all $t\ge0$. In order to obtain a rate of convergence towards zero, we use the recursive inequality
\begin{equation}\label{eq:recursive_E}
E_{t+1} \le E_t -\frac{h\alpha}{J} E_t^2,
\end{equation}
which followed by \eqref{eq:decrease}. We will prove by induction, that $E_t\le \frac{\sigma_u}{(t+1)}$ for $\sigma_u = \frac{J}{\alpha h}$. Firstly, for $t=0$ we have that
$\frac{J}{\alpha h} \ge E_0$ by condition on $h$. Now suppose that the assertion is satisfied for some $t\in\mathbb{N}$, i.e.~$E_t \le \frac{\sigma_u}{(t+1)}$. By \eqref{eq:recursive_E} and the fact that $x\mapsto x-qx^2$ is monotonically increasing for $x\le 1/q$, we obtain with 
$
E_t\le \frac{\sigma_u}{(t+1)} \le \frac{J}{\alpha h},
$
that
\begin{align*}
E_{t+1} &\le \frac{\sigma_u}{t+1} - \frac{\alpha h}{J} \frac{\sigma_u^2}{(t+1)^2}\\
		& = \frac{\sigma_u}{t+2} + \frac{\sigma_u}{(t+1)(t+2)} - \frac{\alpha h}{J}\frac{\sigma_u^2}{(t+1)^2}\\
		&\le \frac{\sigma_u}{t+2} + \frac{\sigma_u - \frac{\alpha h}{J} \sigma_u^2}{(t+1)^2}\\
		&= \frac{\sigma_u}{t+2},
\end{align*}
where we have used $\sigma_u = \frac{J}{\alpha h}$ in the last line. The upper bound on $C_t^{zz}$ follows from the spectral properties of the Frobenius norm.

%Suppose that \sw{$h\le ?$} such that 
%\[C_0^{zz} A_0 C_0^{zz} \succeq hC_0^{zz}A_0 C_0^{zz} A_0 C_0^{zz},\]
%then it follows that $C_0^{zz} \succeq C_1^{zz}$.  
{We pick $h$ small so that $$C_t^{zz} A_{t+1} C_t^{zz} \succeq hC_t^{zz}A_{t+1} C_t^{zz} A_{t+1} C_t^{zz},$$ when $t=0$. Then using induction one can show always holds since $C_t$ decreases.}
%\swtd{Is that true? It depends on $A_t$ right?}

%Hence, the $(C_t^{zz})_{t\ge0}$ solves
%\[C_{t+1}= (I-  hC_t A_{t+1})C_t(I-h C_t A_{t+1})^\top{\leq C_t-h C_t A_{t+1} C_t\leq C_t-\alpha C_t^2} , \quad C_0 = C_0^{zz}\,, \]
%which gives you the upper bound. %Lower bound comes from $C_{t+1}\geq C_t-2h C_t(A_{\max}^2+\alpha) C_t$
%\\
For the lower bound on $C_t^{zz}$, we observe that
\[ C_{t+1}^{zz}\geq C_t^{zz}-2h C_t^{zz}(A_{\max}^2+\alpha) C_t^{zz},\]
since $C_t^{zz}A_{t+1} C_t^{zz} A_{t+1} C_t^{zz}$ is symmetric and positive semi-definite. 
Finally, we apply the upper bound on $\|C_t^{zz}\|_{\mathcal F}\le E_t$ to derive
\begin{align*}
\| C_{t+1}^{zz}-C_t^{zz}\|_{\mathcal F} &\le h^2 (A_{\max}+\alpha)^2 E_t^3 + h(A_{\max}+\alpha)E_t^2\\ &\le ((A_{\max}+\alpha)+h(A_{\max}+\alpha)^2)\frac{J}{\alpha(t+1)^2},
\end{align*}
for sufficiently large $t\ge0$.

\end{proof}
For the remaining analysis, we will denote $C_t = C_t^{zz}$ and suppress the dependency on $z$.
We continue with bounds on the increments $\|\bar z_{t+1} -\bar z_{t}\|^2$ and also uniform moment bounds on the dynamic itself.
\begin{lem}[Moment bound]\label{lem:uniformbounds}
It holds true that there exists $B_1>0$ such that
\[\E[\|\bar z_{t+1} -\bar z_{t}\|^2\mid \mathcal F_t] \le \frac{B_1}{(t+1)^2},\]
and that there exists $B_2>0$ such that
\[\max(\E[\|\bar z_t\|^2],\E[\|\bar z_t - z_\ast\|^2])\le B_2{(1+\log(t))}\,. \]
\end{lem}
\begin{proof}
We have
\begin{align*}
\|\bar z_{t+1} - \bar z_t\|^2 &= h^2 \| C_t (A_{t+1}^{(\alpha)}-S_{t+1}^{\top}S_{t+1} z_\ast + S_{t+1}^\top w_{t+1})\|^2
\\ &\le h^2 \|C_t\|^2 (2\|A_{t+1}^{(\alpha)}-S_{t+1}^{\top}S_{t+1} z_\ast\|^2+2\|S_{t+1}^\top w_{t+1}\|^2),
\end{align*}
and, therefore, by Lemma~\ref{lem:ensemble_collapse}
\begin{align*}
 &\E[\|\bar z_{t+1} - \bar z_t\|^2\mid \mathcal F_t]\\ &\le h^2\frac{\sigma_u^2}{(t+1)^2}\Bigg(2\E[\|A_{t+1}^{(\alpha)}-S_{t+1}^{\top}S_{t+1} z_\ast\|^2\mid \mathcal F_t]+2\E[\|S_{t+1}^\top w_{t+1}\|^2\mid \mathcal F_t]\Bigg)\\ &\le \frac{B_1}{(t+1)^2},
 \end{align*}
for some $B_1>0$ independent of $t$. %\sw{(to prove, the second claim of uniformly bounded second moments. It should directly follow from the first claim.)}
%\\
%\textcolor{blue}{For the second claim, we have
%\begin{align*}
%\max(\E[\|\bar z_t\|^2],\E[\|\bar z_t - z_\ast\|^2]) &= \max(\E[\|\bar z_t\|^2],\E h^2 \| C_t (A_t^{(\alpha)}-S_t^{\top}S_t z_\ast + S_t^\top w_t)\|^2) \\
%& \leq \max(\E[\|\bar z_t\|^2],2h^2\mathbb{E}[ \|C_t\|^2 (\|A_t^{(\alpha)}-S_t^{\top}S_t z_\ast\|^2+\|S_t^\top w_t\|^2)] \\
%&\leq \max(\E[\|\bar z_t\|^2], B_1)
%\\&=: B_2,
%\end{align*}
%which completes the proof.}
{For the second claim, we have
\begin{align*}
\E[\|\bar z_t\|^2] &\le \E[\|\bar z_0\|^2] + t \sum_{k=1}\E[\|\bar z_{k+1}-\bar z_k\|^2]\\ &\le \E[\|\bar z_0\|^2] + t \sum_{k=1}^t \frac{B_1}{(k+1)^2} \\ &\le \E[\|\bar z_0\|^2] + \sum_{k=1}^t \frac{B_1}{(k+1)} \le B_2 (1+\log(t)),
\end{align*}
for some constant $B_2>0$. The bound for  $\E[\|\bar z_t-z_\ast\|^2]$ follows by similar argumentation.  
}
%\swtd{Unfortunately, I believe we obtain another $\log(t)$ factor. At least, I don't see how to avoid it. However, this factor just moves through to the final error bound, and we do not need to adapt the proof below too much.}
\end{proof}

\section{Convergence analysis - proof of the main results}
\label{sec:main}
In this section, we provide our main convergence analysis, which is related to the proofs of Theorems \ref{thm:main} - \ref{thm:main2}.
In particular, each theorem is separated, based on the type of data, that we assume for our problem, where our first proof is based 
on the assumption of using ergodic data, with a corollary following for i.i.d. data, and finally the second theorem assuming periodic data.  We use the notation $a\lesssim b$ denoting $a\le \mathsf{C}b$, $a,b \in \R$, for some constant $\mathsf{C}>0$.
\\
Recall, that $J(z) = \|S(z-z_\ast)\|^2 + \frac{\alpha}2\|z\|^2$ such that $J$ is $\alpha$-strongly convex and $L$-smooth, i.e.~$\nabla J$ is $L$-Lipschitz continuous, where $\alpha$ is the smallest and $L$ is the largest eigenvalue of $A^{(\alpha)}$. %These assumptions are provided through Definitions \ref{def:Ls} - \ref{def:sc}. 
Under smoothness it is well-known that the following descent condition holds
\begin{equation}\label{eq:descentcond}
J(\bar z_{t+1}) \le J(\bar z_t) + \langle \nabla J(\bar z_t) , \bar z_{t+1}-\bar z_t\rangle  + \frac{L}{2}\|\bar z_{t+1} -\bar z_{t}\|^2\,.
\end{equation}
Moreover, using the $\alpha$-strong convexity one can easily derive the Polyak-\L ojasiewicz (PL) inequality 
\[J(z^{\mathrm{opt}}) \geq J(z) - \frac{1}{2\alpha}\| \nabla J((z)\|^2,\quad z\in\R^d\, , \]
where $z^{\mathrm{opt}}$ is the unique minimizer of $J$. 

%\swtd{ToDo: Introduce notation $a\lesssim b$ for $a\le Cb$ for some constant $C>0$. Maybe introduce in the introduction?}

\subsection{Ergodic data convergence proof}

\begin{proof}[Proof of Theorem~\ref{thm:main}]
We proceed from \eqref{eq:descentcond} by considering the following conditional expectation
\begin{align*}
\E[J(\bar z_{t+1})\mid \mathcal F_t] &\le J(\bar z_t) - \langle \nabla J(z_t),h C_t \nabla J(z_t) \rangle\\ &\quad- \langle \nabla J(z_t) , hC_t (A_{t+1}^{(\alpha)}-A^{(\alpha)}) (\bar z_t-z_\ast)\rangle + \frac{L}{2}\E[\|\bar z_{t+1}-\bar z_t\|^2\mid \mathcal F_t]\\
&\le J(\bar z_t) - \langle \nabla J(z_t), hC_t \nabla J(z_t) \rangle\\ &\quad - \langle \nabla J(z_t) , hC_t (A_{t+1}^{(\alpha)}-A^{(\alpha)}) (\bar z_t-z_\ast)\rangle + \frac{L}2 \frac{B_1}{(t+1)^2},
\end{align*}
which holds almost surely. Since $J$ is $\alpha$-strongly convex and $L$-smooth, it satisfies the PL-inequality, which results in
\[\|\nabla J(z)\|^2 \ge q(J(z)-J(z^{\mathrm{opt}})), \]
with $q=2\alpha>0$ and for all $z\in\R^d$. Recall that $A_{t+1}^{(\alpha)}-A^{(\alpha)} = A_{t+1}-A$. We apply the lower bound and upper bounds on $C_t$ to deduce
\begin{align*}
\E[J(\bar z_{t+1})-J(z^{\mathrm{opt}})\mid \mathcal F_t] &\le (1-\frac{\sigma_l}{t+1}q) \left(J(\bar z_t)-J(z^{\mathrm{opt}})\right)\\&\quad - \langle \nabla J(z_t) , hC_t (A_{t+1}-A) (\bar z_t-z_\ast)\rangle + \frac{L}2 \frac{B_1}{(t+1)^2}\\
&=: \rho_t \left(J(\bar z_t)-J(z^{\mathrm{opt}})\right) \\ &\quad- \langle \nabla J(z_t) , hC_t (A_{t+1}-A) (\bar z_t-z_\ast)\rangle + \frac{L}2 \frac{B_1}{(t+1)^2},
\end{align*}
where we have defined $\rho_t := 1-\frac{\sigma_l q}{t+1}$. %Without loss of generality we may assume 
By assumption we have that $\rho_t\in(0,1)$, using $0<\lambda = \sigma_l q<1$ for sufficiently small $\sigma_l$. Taking expectation and defining 
$$
\Delta_t = J(\bar z_{t})-J(z^{\mathrm{opt}}),\quad b_t = -\langle \nabla J(z_t) , hC_t (A_{t+1}-A) (\bar z_t-z_\ast)\rangle,
$$
we obtain the iterative bound 
\begin{align*}
\E[\Delta_T] \le \left(\prod_{s=0}^{T-1} \rho_s\right)\E[\Delta_0] + \sum_{s=0}^{T-1} \frac{\prod_{r=s}^{T-1}\rho_{r+1}}{\rho_T} \E[b_s] + \sum_{s=0}^{T-1} \frac{\prod_{r=s}^{T-1}\rho_{r+1}}{\rho_T} \frac{L}2 \frac{B_1}{(s+1)^2}.
\end{align*}
By discrete Gronwall's inequality we have that
\begin{align}
\label{eq:grom_deltat}
\E[\Delta_T] \le \left(\frac{1}{T}\right)^{\lambda} \E[\Delta_0] + \sum_{s=0}^{T-1} \left(\frac{s+1}{T}\right)^{\lambda} \E[b_s] + \sum_{s=0}^T \left(\frac{s+1}{T}\right)^{\lambda} \frac{L}2 \frac{B_1}{(s+1)^2}\,.
\end{align}
%\swtd{One last problem I have recognized: Somehow we were beating the convergence rate of SGD in the strongly convex setting for $\lambda>1$. Since this is not possible (there are counterexamples), I believe we need $\lambda\le1$ for the Gronwall argument, since $\rho_t\in(0,1)$ needs to be satisfied. I have included this as assumption.}
%where \sw{$\lambda = \frac{1}{2}\sigma_l q\le1$ for sufficiently small $\sigma_l$}. 
The second term of the upper bound writes as
\begin{align*}
\sum_{t=0}^{T-1} \left(\frac{s+1}{T}\right)^{\lambda} \E[b_s]  &= \E\Bigg[\sum_{t=0}^{T-1} \left(\frac{t+1}{T}\right)^{\lambda} \Big\langle A^{(\alpha)} \bar z_t - Az_\ast, C_t (A-A_{t+1})(\bar z_t-z_\ast)\Big\rangle \Bigg]\\
&= \E\Bigg[\sum_{t=0}^{T-1} \left(\frac{t+1}{T}\right)^{\lambda} (\bar z_t - z_\ast)^\top A^\top C_t (A-A_{t+1})(\bar z_t-z_\ast) \Bigg]\\
&+\E\Bigg[\sum_{t=0}^{T-1} \left(\frac{t+1}{T}\right)^{\lambda} \alpha\bar z_t C_t (A-A_{t+1})(\bar z_t-z_\ast) \Bigg].
\end{align*}
We will consider only the first sum, bounding the second sum will follow similarly. Firstly, observe that with $g_t := \bar z_t-z_\ast$ we can write
\begin{align*}
\E\Bigg[\sum_{t=0}^{T-1} \left(\frac{t+1}{T}\right)^{\lambda} &g_t^\top A^\top C_t (A-A_{t+1})g_t \Bigg] \\ &= \E\Bigg[\sum_{t=0}^{T-1} \left(\frac{t+1}{T}\right)^{\lambda} g_t^\top A^\top C_t (A-A_{t+1+\tau_\varepsilon})g_t \Bigg]\\
&+\E\Bigg[\sum_{t=0}^{T-1} \left(\frac{t+1}{T}\right)^{\lambda} g_t^\top A^\top C_t (A_{t+1+\tau_\varepsilon}-A_{t+1})g_t \Bigg]\,.
\end{align*}
Since $g_t$ and $C_t$ are $\mathcal F_t$-measurable it follows that
\begin{align*}
\E\Bigg[\sum_{t=0}^{T-1} \left(\frac{t+1}{T}\right)^{\lambda} &g_t^\top A^\top C_t (A-A_{t+1+\tau_\varepsilon})g_t \Bigg]\\ &= \sum_{t=0}^{T-1} \left(\frac{t+1}{T}\right)^{\lambda} \E\Bigg[g_t^\top A^\top C_t \E[(A-A_{t+1+\tau_\varepsilon})\mid \mathcal F_t] g_t \Bigg]\\
&\le \sum_{t=0}^{T-1} \left(\frac{t+1}{T}\right)^{\lambda} \E\Bigg[\|g_t\|^2 \|A\| \|C_t\| \|A-\E[A_{t+1+\tau_\varepsilon}\mid \mathcal F_t]\| \Bigg]\\
&\le \varepsilon \sum_{t=0}^{T-1} \left(\frac{t+1}{T}\right)^{\lambda} \Bigg(\frac{\sigma_u}{\alpha(t+1)} \Bigg)\|A\| \E[\|g_t\|^2]\\
&{\le \varepsilon \sum_{t=0}^{T-1} \left(\frac{t+1}{T}\right)^{\lambda} \Bigg(\frac{\sigma_u}{\alpha(t+1)} \Bigg)\|A\| B_2 (1+\log(t))}\\
&\lesssim \varepsilon ,
\end{align*}
where %$a \lesssim b \implies a \leq Cb$, and that 
we have used Lemma~\ref{lem:ensemble_collapse} and the assumption that $$\|A-\E[A_{t+1+\tau_\varepsilon}\mid \mathcal F_t]\|\le \varepsilon.$$ 
%\swtd{ToDo: Introduce notation $a\lesssim b$ for $a\le Cb$ for some constant $C>0$. Maybe introduce in the introduction?}
We continue with the bound
\begin{align*}
\E\Bigg[\sum_{t=0}^{T-1} \left(\frac{t+1}{T}\right)^{\lambda} &g_t^\top A^\top C_t (A_{t+1+\tau_\varepsilon}-A_{t+1})g_t \Bigg]\\ &= {\mathrm{residual}}(\tau_\varepsilon)\\ & + \E\Bigg[\sum_{t=\tau_\varepsilon+1}^T\left(\frac{t+1-\tau_\varepsilon}{T}\right)^{\lambda} g_{t-\tau_\varepsilon}^\top A^\top C_{t-\tau_\varepsilon} A_{t+1} g_{t-\tau_\varepsilon} \Bigg]\\ & - \E\Bigg[\sum_{t=\tau_\varepsilon+1}^T\left(\frac{t+1}{T}\right)^{\lambda} g_t^\top A^\top C_t A_{t+1} g_t \Bigg],
\end{align*}
where we define the residual to be
\begin{equation}\label{eq:residuals}
\begin{split}
{\mathrm{residual}}(\tau_\varepsilon)&:= \E\Bigg[\sum_{t=T-1+1}^{T-1+\tau_\varepsilon} \left(\frac{t+1-\tau_\varepsilon}{T}\right)^\lambda g_{t-\tau_\varepsilon}^\top A^\top C_{t-\tau_\varepsilon} A_{t+1} g_{t-\tau_\varepsilon}\Bigg] \\&- \E\Bigg[\sum_{t=0}^{\tau_\varepsilon} \left(\frac{t+1}{T}\right)^\lambda g_t^{\top} A^\top C_t A_{t+1} g_t\Bigg], 
\end{split}
\end{equation}
which we discuss in the later part of the proof. We apply the estimate
\begin{align*}
&\E\Bigg[\Bigg|\left(\frac{t+1-\tau_\varepsilon}{T}\right)^{\lambda} g_{t-\tau_\varepsilon}^\top A^\top C_{t-\tau_\varepsilon} A_{t+1} g_{t-\tau_\varepsilon} ]  - \left(\frac{t+1-\tau_\varepsilon}{T}\right)^{\lambda} g_t^\top A^\top C_t A_{t+1} g_t\Bigg|\Bigg]\\ &\le \left(\frac{t+1}{T}\right)^\lambda \|A_{t+1}\| \|A^{(\alpha)}\| \Bigg(\E[\|g_{t-\tau_\varepsilon}-g_t\| \|g_t\|\|C_t\|]  + \E[\|g_t\|^2\|C_t-C_{t-\tau_\varepsilon}\|]  \\& \quad+ \E[\|g_{t-\tau_\varepsilon}\| \|g_{t-\tau_\varepsilon}-g_t\|\|C_t\|] \Bigg)\\
&\lesssim \left(\frac{t+1-\tau_\varepsilon}{T}\right)^\lambda \frac{{(1+\log(t))\tau_\varepsilon}}{(t+1)(t+1-\tau_\varepsilon)}\,,
\end{align*}
where we have used Hölder's inequality to uniformly bound $\E[\|g_{t-\tau_\varepsilon}\| \|g_{t-\tau_\varepsilon}-g_t\|]$ and $\E[\|g_{t-\tau_\varepsilon}-g_t\| \|g_t\|]$ by Lemma~\ref{lem:uniformbounds}. {More precisely, we have used that 
$$
\E[\|g_{t-\tau_\varepsilon}-g_t\|^2]^{1/2}\le \sum_{k=t-\tau_\varepsilon}{t} \E[\|\bar z_k-\bar z_{k-1}\|^2 \lesssim \frac{\tau_\varepsilon}{t+1-\tau_\varepsilon}.
$$}
Hence, it follows that
\begin{align*}
\Bigg|\E \Bigg[\sum_{t=\tau_\varepsilon+1}^T & \left(\frac{t+1-\tau_\varepsilon}{T}\right)^{\lambda} g_{t-\tau_\varepsilon}^\top A^\top C_{t-\tau_\varepsilon} A_{t+1} g_{t-\tau_\varepsilon} \Bigg] \\ & - \E\Bigg[\sum_{t=\tau_\varepsilon+1}^T\left(\frac{t+1}{T}\right)^{\lambda} g_t^\top A^\top C_t A_{t+1} g_t \Bigg] \Bigg|\\ & \lesssim {\tau_\varepsilon} \sum_{t=0}^T \left(\frac{t+1}{T}\right)^\lambda \frac{{1+\log(t)}}{(t+1)^2} = \frac{{\tau_\varepsilon}}{(T+1)^\lambda} \sum_{t=0}^T \frac{{1+\log(t)}}{(t+1)^{2-\lambda}}\,.
\end{align*}
%\swtd{Here we have to consider different settings for $\lambda$. E.g.~if $\lambda=1$, we have  $$\sum_{t=0}^T \frac{1}{(t+1)^{2-\lambda}}\le (1+\log(T)).$$}
For the residual terms in $ {\mathrm{residual}}(\tau_\varepsilon)$ we can use the bound
\begin{align*}
\Bigg| \E\Bigg[\sum_{t=T-1+1}^{T-1+\tau_\varepsilon} &\left(\frac{t+1-\tau_\varepsilon}{T}\right)^\lambda g_{t-\tau_\varepsilon}^\top A^\top C_{t-\tau_\varepsilon} A_{t+1} g_{t-\tau_\varepsilon}\Bigg]\Bigg|\\ &\le a_1 \|A\| {(1+\log(T-1))B_2} \sum_{t=T-\tau_\varepsilon}^{T-1} \left(\frac{t+1}{T}\right)^\lambda\|C_t\|,
\end{align*}
with 
\[a_1 = \max_{t=T,\dots,T-1+\tau_\varepsilon}\ \|A_{t+1}\|,\quad \text{and}\quad \max_{T-\tau_\varepsilon,\dots, T-1}\ \E[\|g_t\|]{\le B_2(1+\log(T-1))} \,, \]
and similarly
\begin{align*}
\Bigg|\E\Bigg[\sum_{t=0}^{\tau_\varepsilon} \left(\frac{t+1}{T}\right)^\lambda g_t^{\top} A^\top C_t A_{t+1} g_t\Bigg]\Bigg| \le b_1 \|A\| {(1+\log(\tau_\varepsilon))B_2} \sum_{t=0}^{\tau_\varepsilon} \left(\frac{t+1}{T}\right)^\lambda \|C_t\|,
\end{align*}
with 
\[ b_1 = \max_{t=0,\dots,\tau_\varepsilon}\ \|A_{t+1}\|,\quad  \text{and}\quad \max_{t=0,\dots,\tau_\varepsilon}\ \E[\|g_t\|]{\le B_2(1+\log(\tau_\varepsilon))}\,.\]
This means, that the residual has the following bound% ${\mathrm{residual}}(\tau)\lesssim \left(\frac{1}{T}\right)^{\lambda}$ 
\begin{align*}
\E\Bigg[\sum_{t=T-1+1}^{T-1+\tau_\varepsilon} &\left(\frac{t+1-\tau_\varepsilon}{T}\right)^\lambda g_{t-\tau_\varepsilon}^\top A^\top C_{t-\tau_\varepsilon} A_{t+1} g_{t-\tau_\varepsilon}\Bigg]\\ &- \E\Bigg[\sum_{t=0}^{\tau_\varepsilon} \left(\frac{t+1}{T}\right)^\lambda g_t^{\top} A^\top C_t A_{t+1} g_t\Bigg]\\ &\lesssim \left(\frac{1}{T}\right)^{\lambda} \tau_\varepsilon (1+\log(\tau_\varepsilon)+\log(T))\,.
\end{align*}
Finally, using that
\begin{equation}\label{eq:finitesum}
\sum_{t=0}^T \left(\frac{t+1}{T}\right)^\lambda \frac{1}{(t+1)^2} \lesssim \frac{1}{T^\lambda},
\end{equation}
for $\lambda<1$, we obtain a total error bound of
\[ \E[\Delta_t] \lesssim \frac{(1+\tau_\varepsilon (1+\log(\tau_\varepsilon)+\log(T)))}{T^\lambda} %\left(\frac{1}{T}\right)^{\lambda}{(1+\tau_\varepsilon (1+\log(\tau_\varepsilon)+\log(T)))} + {\tau_\varepsilon (1+\log(T))} \sum_{t=0}^T \left(\frac{t+1}{T}\right)^\lambda \frac{1}{(t+1)^2} 
+ \varepsilon \,.\]
\end{proof}
%\begin{remark}
%If we consider the special case of $\lambda=1$, then we have the following bound $$\sum_{t=0}^T \frac{1}{(t+1)^{2-\lambda}}\le (1+\log(T)).$$
%\end{remark}
The proof of Corollary~\ref{thm:main3} even simplifies crucially, since we have that $A_t=S_t^\top S_t$ are independent and unbiased estimators of $A$.
\begin{proof}[Proof of Corollary~\ref{thm:main3}]
The assertion directly implies from the fact that $\E[S_{t+1}^\top S_{t+1}\mid \mathcal F_t] = \E[A_{t+1}\mid \mathcal F_t] = A$ for any $t\ge0$, which simplifies equation \eqref{eq:grom_deltat} to
\begin{equation*}
\E[\Delta_T] \le \left(\frac{1}{T}\right)^{\lambda} \E[\Delta_0] + \frac{B_1L}2\sum_{s=0}^T \left(\frac{s+1}{T}\right)^{\lambda}  \frac{1}{(s+1)^2},
\end{equation*}
with $\lambda<1$. The assertion follows again using \eqref{eq:finitesum}.
%\[ \sum_{t=0}^T \left(\frac{t+1}{T}\right)^\lambda \frac{1}{(t+1)^2} \lesssim \frac{1}{T^\lambda}\,.\]
\end{proof}
%\begin{remark}
%\textcolor{blue}{Maybe just briefly discuss this setting, without proving anything? To discuss.}
%\swtd{Here we have to consider different settings for $\lambda$. E.g.~if $\lambda=1$, we have  $$\sum_{t=0}^T \frac{1}{(t+1)^{2-\lambda}}\le (1+\log(T)).$$}
%\swtd{This is no extension. We have to control the sum 
%$$
%\sum_{s=0}^T \left(\frac{s+1}{T}\right)^2 \frac{1}{(s+1)^2},
%$$ in general for $\lambda>0$. However, $\lambda$ can be arbitrarily small.}
%\end{remark}
%
%\subsection{i.i.d. data convergence proof}
%
%\swtd{ToDo: formulate Theorem~\ref{thm:main3} as corollary of Theorem~\ref{thm:main}. We then no longer need a proof of this statement.}
%\begin{proof}[Proof of Theorem \ref{thm:main3}]
%\textcolor{blue}{TBD}
%\end{proof}
\subsection{Periodic data convergence proof}
%\subsection{Old notes}
\begin{proof}[Proof of Theorem~\ref{thm:main2}]
 %We now proceed in a similar fashion as above, which
 %is given through the following lemma.

%The extension should also be easy, everything remains same till here
%\begin{align*}
%\sum_{n=1}^N \Big(\frac{n}{N}\Big)^{\lambda}v_n^TAC_n (A-A_n) v_n
%&=\sum_{n=1}^N \Big(\frac{n}{N}\Big)^{\lambda}v_n^TAC_n (A-A_{n+\tau}) v_n\\&+
%\sum_{n=1}^N \Big(\frac{n}{N}\Big)^{\lambda}v_n^TAC_n (A_{n+\tau}-A_{n}) v_n.
%\end{align*}
%We replace it with 
%\begin{align*}
%\sum_{n=1}^N \Big(\frac{n}{N}\Big)^{\lambda}v_n^TAC_n (A-A_n) v_n
%&=\sum_{n=1}^N \Big(\frac{n}{N}\Big)^{\lambda}v_n^TAC_n (A-\frac{1}{T}A_{n:n+T}) v_n\\&+
%\frac1T\sum_{\tau=1}^T\sum_{n=1}^N \Big(\frac{n}{N}\Big)^{\lambda}v_n^TAC_n (A_{n+\tau}-A_{n}) v_n.
%\end{align*}
%The first part, when we do conditional expectation. The second part can be bound using the methods above for each individual $\tau=1,...,T$.

The proof follows the similar lines as the proof of Theorem~\ref{thm:main}. Proceeding in a similar fashion, 
everything remains the same up until equation \eqref{eq:grom_deltat}. Recall that we have 
\begin{align*}
\sum_{t=0}^{T-1} \left(\frac{s+1}{T}\right)^{\lambda} \E[b_s]  &= \E\Bigg[\sum_{t=0}^{T-1} \left(\frac{t+1}{T}\right)^{\lambda} \Big\langle A^{(\alpha)} \bar z_t - Az_\ast, C_t (A-A_{t+1})(\bar z_t-z_\ast)\Big\rangle \Bigg]\\
&= \E\Bigg[\sum_{t=0}^{T-1} \left(\frac{t+1}{T}\right)^{\lambda} (\bar z_t - z_\ast)^\top A^\top C_t (A-A_{t+1})(\bar z_t-z_\ast) \Bigg]\\
&+\E\Bigg[\sum_{t=0}^{T-1} \left(\frac{t+1}{T}\right)^{\lambda} \alpha\bar z_t C_t (A-A_{t+1})\bar z_t-z_\ast) \Bigg],
\end{align*}
where we again only consider the first term. The second term will result in similar bounds. In the periodic case we now proceed differently to the proof of Theorem~\ref{thm:main}, and decompose the expression in the following way
\begin{align*}
\E\Bigg[\sum_{t=0}^{T-1} \left(\frac{t+1}{T}\right)^{\lambda} &g_t^\top A^\top C_t (A-A_t)g_t \Bigg]\\ &= \E\Bigg[\sum_{t=0}^{T-1} \left(\frac{t+1}{T}\right)^{\lambda} g_t^\top A^\top C_t \Big(A-\frac{1}{\tau_{\varepsilon}}A_{(t+1):(t+1+\tau_{\varepsilon})}\Big)g_t \Bigg]\\
&+\E\Bigg[\frac{1}{\tau_{\varepsilon}} \sum^{\tau_{\varepsilon}-1}_{k=0} \sum_{t=0}^{T-1} \left(\frac{t+1}{T}\right)^{\lambda} g_t^\top A^\top C_t (A_{t+1+k}-A_t)g_t \Bigg].
\end{align*}
As before, since $g_t$ and $C_t$ are $\mathcal F_t$-measurable it follows that
\begin{align*}
\E\Bigg[\sum_{t=0}^{T-1}& \left(\frac{t+1}{T}\right)^{\lambda} g_t^\top A^\top C_t \Big(A-\frac{1}{\tau_{\varepsilon}}A_{(t+1):(t+1+\tau_{\varepsilon})}\Big)g_t \Bigg] \\&= \sum_{t=0}^{T-1} \left(\frac{t+1}{T}\right)^{\lambda} \E\Bigg[g_t^\top A^\top C_t \E[(A-\frac{1}{\tau_{\varepsilon}}A_{(t+1):(t+1+\tau_{\varepsilon})}\mid \mathcal F_t] g_t \Bigg]\\
&\le \sum_{t=0}^{T-1} \left(\frac{t+1}{T}\right)^{\lambda} \E\Bigg[\|g_t\|^2 \|A\| \|C_t\| \|A-\E[\frac{1}{\tau_{\varepsilon}}A_{(t+1):(t+1+\tau_{\varepsilon})}\mid \mathcal F_t]\| \Bigg]\\
&\le \varepsilon \sum_{t=0}^{T-1} \left(\frac{t+1}{T}\right)^{\lambda} \Bigg(\frac{\sigma_u}{\alpha(t+1)} \Bigg)\|A\| \E[\|g_t\|^2]\\
&{\le \varepsilon \sum_{t=0}^{T-1} \left(\frac{t+1}{T}\right)^{\lambda} \Bigg(\frac{\sigma_u}{\alpha(t+1)} \Bigg)\|A\| B_2 (1+\log(t))}\\
&\lesssim \varepsilon,
\end{align*}
where we have made use of Assumption \ref{assum:period} and Lemma \ref{lem:ensemble_collapse}.
Now we consider
\begin{align*}
\E\Bigg[\frac{1}{\tau_{\varepsilon}} &\sum^{\tau_{\varepsilon}-1}_{k=0} \sum_{t=0}^{T-1} \left(\frac{t+1}{T}\right)^{\lambda} g_t^\top A^\top C_t (A_{t+1+k}-A_t)g_t \Bigg] \\ &= \frac{1}{\tau_{\varepsilon}}\sum_{k=0}^{\tau_{\varepsilon}-1}\mathrm{residual}(k) + \E\Bigg[\frac{1}{{\tau_{\varepsilon}}}\sum^{\tau_{\varepsilon}-1}_{k=0}\sum_{t=k+1}^T\left(\frac{t+1-k}{T}\right)^{\lambda} g_{t-k}^\top A^\top C_{t-k} A_{t+1} g_{t-k} \Bigg]\\ & - \E\Bigg[\frac{1}{{\tau_{\varepsilon}}}\sum^{\tau_{\varepsilon}-1}_{k=0}\sum_{t=k+1}^T\left(\frac{t+1}{T}\right)^{\lambda} g_t^\top A^\top C_t A_{t+1} g_t \Bigg],
\end{align*}
where we again use the residual defined in \eqref{eq:residuals}. With similar arguments as in the proof of Theorem~\ref{thm:main} replacing each bound applied to fixed $\tau_\varepsilon$ now by $k=0,\dots,\tau_{\varepsilon}-1$, we obtain
\begin{align*}
\Bigg|\E\Bigg[\frac{1}{\tau_{\varepsilon}} \sum^{\tau_{\varepsilon}-1}_{k=0} \sum_{t=0}^{T-1} &\left(\frac{t+1}{T}\right)^{\lambda} g_t^\top A^\top C_t (A_{t+1+k}-A_t)g_t \Bigg]\Bigg|\\&\lesssim \frac{1}{\tau_{\varepsilon}}\sum_{k=0}^{\tau_{\varepsilon}-1} k \sum_{t=0}^T \left(\frac{t+1}{T}\right)^\lambda \frac{1+\log(t)}{(t+1)^{2}}
\\ &\quad+ \frac{1}{T^\lambda} \frac{1}{\tau_{\varepsilon}}\sum_{k=0}^{\tau_\varepsilon-1} k\big(1+\log(k)\big)+\log(T)\\
&\lesssim   (1+\log(\tau_{\varepsilon}))(1+\log(T)) \sum_{t=0}^T \left(\frac{t+1}{T}\right)^\lambda \frac{1}{(t+1)^2} 
\\&+ \frac{\log(T)}{T^\lambda}+\frac{(1+\log(\tau_{\varepsilon}))^2}{T^\lambda}\,.
\end{align*}
Again, with \eqref{eq:finitesum} the final error bound is then given by
\begin{align*}
 \E[\Delta_t] &\lesssim \left(\frac{1}{T}\right)^{\lambda} (1+(1+\log(\tau_\varepsilon))^2+\log(T))\\&\quad + (1+\log(\tau_\varepsilon)) (1+\log(T)) \sum_{t=0}^T \left(\frac{t+1}{T}\right)^\lambda \frac{1}{(t+1)^2} + \varepsilon\\
 &\lesssim \frac{1+(1+\log(\tau_\varepsilon))^2+(1+\log(\tau_\varepsilon))(1+\log(T))}{T^\lambda}\,.
 \end{align*}
%where we define the residual to be
%\begin{align*}
%{\mathrm{residual}}(\tau)&:= \E\Bigg[\frac{1}{{\tau_{\varepsilon}}}\sum ^{\tau_{\varepsilon}}_{k=0}\sum_{t=T-1+1}^{T-1+\tau} \left(\frac{t+1-\tau}{T}\right)^\lambda g_{t-\tau}^\top A^\top C_{t-\tau} A_t g_{t-\tau}\Bigg] \\&- \E\Bigg[\frac{1}{{\tau_{\varepsilon}}}\sum ^{\tau_{\varepsilon}}_{k=0}\sum_{t=0}^\tau \left(\frac{t+1}{T}\right)^\lambda g_t^{\top} A^\top C_t A_t g_t\Bigg], 
%\end{align*}
%\nc{where now we aim to bound using the methods above for each individual $\tau=1,...,T$.
%}
\end{proof}

\section{Numerical experiments}
\label{sec:num}

In this section, we introduce numerical toy models for which we aim to verify our findings from the previous sections, and how
our methodology performs in general. Specifically, we will present and implement two different numerical experiments based on Example~\ref{ex:darcyflow}. We consider two particular setups, the first where we split up our domain into sub-domains to construct a dynamic observation operator, and the second where we consider a dynamic solution operator. Our simulations will be based on the 2D Darcy flow Poisson equation. \textcolor{black}{We will also compare the effect of the ensemble size on the convergence}.

%\subsection{Darcy flow with discrete observation points}
Recall, the model of interest considered in Example~\ref{ex:darcyflow} is given as $2$-dimensional elliptic PDE model
\begin{equation*}%\label{eq:darcy_flow}
\begin{cases}
-\nabla \cdot (\exp(a)\nabla p) = z^\ast,& x\in D\\
p=0,& x\in\partial D
\end{cases}\, ,
\end{equation*}
with domain $D:=(0,1)^2$ and subject to zero Dirichlet boundary conditions. Our aim is to recover the unknown source term $z^\ast\in L^\infty(D)$ from different observation models using discrete observation points of the solution $p\in \mathcal V:= H_0^1(D)\cap H^2(D)$. We will test the different models presented in Example~\ref{ex:darcyflow}. Recall, that $\mathcal O_{[x_1:x_K]}:\mathcal V\to \R^K$ evaluates the solution $p\in\mathcal V$ in $K$ randomly picked observation points and 
given $a\in L^\infty$ the linear operator $G_a:L^\infty(D)\to \mathcal V$ solves the equation \eqref{eq:darcy_flow}.
The dynamic forward model is then given by $S_t = \mathcal O_t \circ G_t$, where we specify the different choices for $\mathcal O_t$ and $G_t$ in the following. Our PDE is solved using a centred finite difference method with specified mesh size of $h_*=1/100$, where we have that $d=100$.
%\swtd{conflict of notation, $h$ is already defined for the scaling of the steps in the dynamic EKI algorithm.}
%discrete observation points of the solution $p\in \mathcal V:= H_0^1(D)\cap H^2(D)$. Let $G:L^\infty(D)\to \mathcal V$ be the solution operator of \eqref{eq:darcy_flow} and define the observation operator to be an linear operator $\mathcal O_{[x_1:x_K]}:\mathcal V\to \R^K$ that evaluates $p\in\mathcal V$ in $K$ randomly picked observation points. More precisely, this means $\mathcal O_{[x_1:x_K]} p = (p(x_1),\dots, p(x_K))^\top\in\R^K$.

\subsection{Periodic experiment: Dynamic observation operator} 
Our first numerical example is a scenario of a dynamic observation operator $\mathcal O_t$ and static solution operator $G$. We consider two different dynamic observation models, where the first one is simply an i.i.d. observation model, the second one takes periodic measurements on sub-domains. %and the third one also picks measurements on sub-domains but with ergodic behavior. 
For both models we assume that the diffusion coefficient $a\in L^\infty(D)$ is fixed. % is to consider to case where we compare the use of periodic data, ergodic data where our construction of the perioidic data is based on splitting up the domain. 
\subsubsection{Independent and identically distributed(i.i.d.) observation model}
We start by defining our observation model to be %decomposed into 
\[ u_t = S_t z^\ast + w_t\,,\]
where $S_t = \mathcal O_t\circ G_a$ and $\mathcal{O}_t = \mathcal{O}_{[x_1(t):x_K(t)]}$ with $x_1(t),\dots,x_K(t)\sim \mathcal U(D)$ are drawn independently with uniform distribution over the entire domain $D$.  %In this case $\mathcal O_t\circ G$ satisfies Assumption \eqref{assum:iid}.

We will use this model to construct our reference solution by defining the empirical objective function for a large number $N\gg 1$ of i.i.d. observations. The empirical objective function is defined by 
\[ z_{\mathrm{ref}} = \arg\min_{z}\ J_N(z),\quad J_N(z) := \frac{1}{2N}\sum_{t=1}^N \|S_t z - u_t\|^2 + \frac{\alpha}2\|z\|^2\, .\]
For our experiments we fix a regularization parameter of $\alpha=2$.

%\nc{Also simon related to the imaging comment, if Xin believes the nonlinear extension is doable, then this could potentially be applied to some imaging problems making it stronger both for theory and numerics. For future work, I mean.}

\subsubsection{Periodic observation model}
In the next example, we consider the periodic observation model from Example~\ref{ex:darcyflow}, (ii), where we decomposed the domain $D$ into $p$ disjoint subsets $D_1,\dots,D_p$ such that $D=\cup_{i=1}^p D_i$. The dynamic observation operator is $\mathcal O_t$ is defined in the following way. For $t=pk+i$ with $k\ge0$ and $i=1,\dots,p$ we draw $x_1(t),\dots, x_K(t)\sim\mathcal U(D_i)$ independently with uniform distribution over the sub-domain $D_i\subset D$. The observation model then reads as
\[ u_t = S_t z^\ast + w_t\,,\]
where $S_t = \mathcal O_t\circ G_a$.

\begin{figure}[h!]
\centering
\hspace*{-1cm}                                          
\includegraphics[scale=0.35]{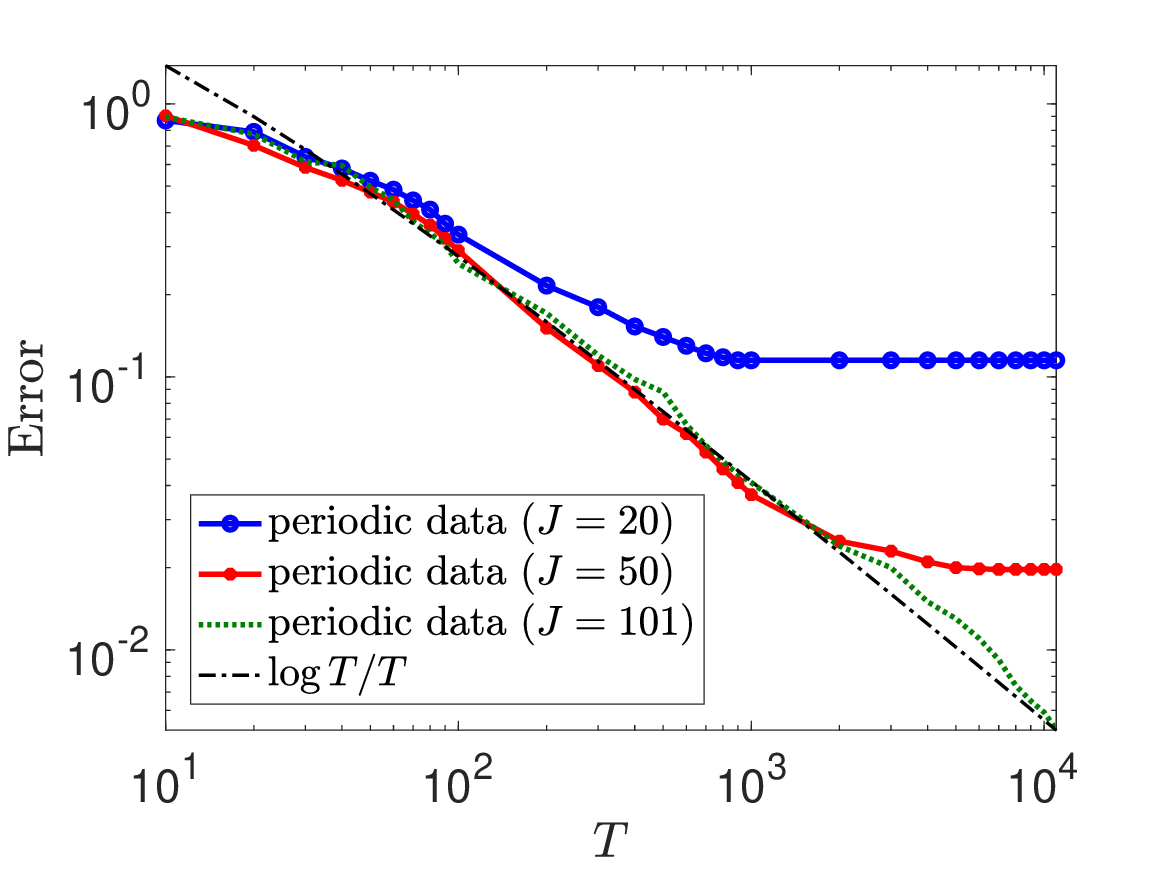}~~\includegraphics[scale=0.35]{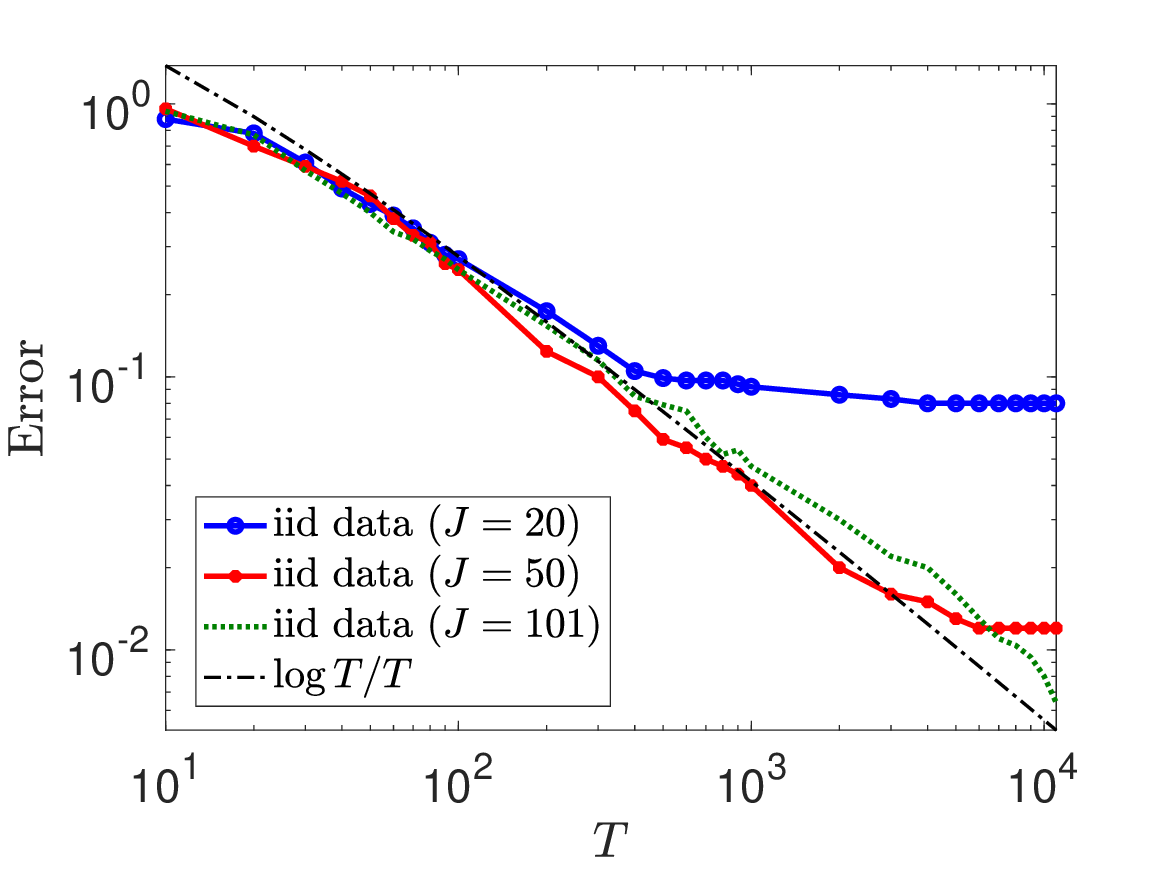}
\caption{Convergence plot for EKI using periodic (left) and i.i.d.~(right) data. \textcolor{black}{We plot the error w.r.t.~the reference solution $z_{\mathrm{ref}}$ vs time $T$.}}
\label{fig:cov}
\end{figure}

Our \textcolor{black}{ground truth} $z^\ast$ will be based on a the Karhunen-Lo\`{e}ve expansion (KLE) which is used to simulate Gaussian random fields. Specifically we consider $z^\ast$ as realization of a Gaussian unknown $\mu\sim \mathcal{N}(0,\mathcal{C})$ with Mat\'{e}rn covariance function, 
\begin{equation}
\label{eq:matern}
\mathcal{C}(x,x') = \frac{2^{1-\nu}}{\Gamma(\nu)} K_{\nu}\bigg(\frac{|x-x'|}{\ell}\bigg)\bigg(\frac{|x-x'|}{\ell}\bigg)^{\nu},
\end{equation}
%\red{[Please check if correct. Look strange]}
and KLE defined as
\begin{equation}
\label{eq:kle}
\mu(x) = \sum^{\mathcal{J}}_{j=1} \phi_j \sqrt{\lambda_j} \xi_j, \quad \xi_j \sim \mathcal{N}(0,I),
\end{equation}
%\red{$J$ collides with ensemble size.}
where $(\phi_j,\lambda_j)$ is the corresponding eigenbasis of $\mathcal{C}$, and $(\nu, \ell) \in \R^+ \times \R$ are
associated hyperparameters, $K_{\nu}(\cdot)$ represents a Bessel function of the second kind and $\Gamma(\cdot)$ is a Gamma function. Specifically our \textcolor{black}{ground truth} $z^{\ast}$ will be chosen with hyperparameter choices of $(\ell,\nu) = (0.07,3.4)$. 
For our periodic data, we split the domain into 10 domains $D_1,\ldots, D_{10}$.  When running our DEKI algorithm we specify $J=50$ ensemble members where we place $T=10,000$. We provide a convergence plot which is given in Figure~\ref{fig:cov}.
\par
\textcolor{black}{From the numerical subplots in Figure~\ref{fig:cov} we see that we attain the theoretical rates of Corollary~\ref{thm:main3} and Theorem~\ref{thm:main2}, for both  the i.i.d. and periodic data in the regime of $J>d$, when $J=101$. However, as expected due to the subspace property Lemma~\ref{lem:subspace}, when we consider the regime $J<d$, where we exemplary choose $J=\{20,50\}$, we notice for both sets of data that the error reaches a plateau. %which implies we do not see convergence towards zero. 
This behavior, as well as the observed improvement in increasing $J$, can be explained by the findings of \cite{CST19, W2022}, which demonstrate that the EKI can be viewed as an optimizer constrained to the subspace spanned by the initial ensemble.}

%\begin{figure}[h!]
%\centering
%\includegraphics[scale=0.35]{iid.eps}
%\caption{Convergence plot for EKI using periodic data. \textcolor{blue}{We plot the error w.r.t.~the reference solution $z_{\mathrm{ref}}$ vs time $T$.}}
%\label{fig:cov}
%\end{figure} 

\subsection{Ergodic experiment: Dynamic solution operator}
Our second experiment uses a static observation operator $\mathcal O$ but incorporates a dynamic solution operator $G_t$ for \eqref{eq:darcy_flow} using different realizations of the diffusion coefficient $a$. \textcolor{black}{For this experiment we solely consider an ergodic observation models, however we note that we can also consider an i.i.d. observation model}. %Since the application is motivated by the ergodic observation model, we will start with the description of the ergodic model first.

We assume that we are in the scenario of Example~\ref{ex:darcyflow}, (iii), where the observation operator is fixed $\mathcal O = \mathcal O_{[x_1:x_K]}$ for fixed locations of observation $x_1,\dots,x_K\in D$ and the diffusion coefficient is unknown. To be more precise, we assume that in general we have access to some statistical information about $a$ in form of a probability distribution $\pi$. For example, this information may come from a pre-stage Bayesian experiment. In many practical scenarios the explicit computation of $\pi$ or the generation of exact samples of $\pi$ is infeasible and for our ergodic observation model we assume that the information about the diffusion coefficient has been generated by a Markov chain Monte Carlo (MCMC) algorithm. As comparison we also consider the simplified setting, where one can generate i.i.d. samples of $\pi$.

%\subsubsection{Independent and identically distributed observation model}

%Firstly, we assume that there is access to a sample of i.i.d. realizations of $\pi$ denoted by $(a_t^{\mathrm{i.i.d.}})_{t\ge0}$. The observations are then defined by 
%\[u_t = S_t z^\ast + w_t\,, \]
%where $S_t = \mathcal O \circ G_t$ with $G_t = G_{a_t^{\mathrm{i.i.d.}}}$. Again, we make use of this model to construct our reference solution
%\[z_{\mathrm{ref}} := \arg\min_{z}\ J_N(z),\quad J_N(z):= \frac{1}{2N}\sum_{t=1}^N \|S_t z - u_t\|^2 + \frac{\alpha}2\|z\|^2\]  
%for a large number of $N\gg 1$ i.i.d. observations. Similarly as before, we fix a regularization parameter of $\alpha = 2$.

%\subsubsection{Periodic observation model}
%Finally, to also compare to the periodic observation model in this scenario, we introduce a periodic sequence of observations through the empirical objective $J_N$ which we have used to construct our reference solution. We split the data set $\{u_t,t=1,\dots,N\}$ generated by the iid observation model into $K$ sub-blocks
%\swtd{stopped here}

\subsubsection{Ergodic observation model}
%We assume that we are in scenario of Example~\ref{ex:darcyflow}, (iii), where the observation operator is fixed $\mathcal O = \mathcal O_{[x_1:x_K]}$ for fixed observation points $x_1,\dots,x_K\in D$ and
As discussed above we assume that the information about the diffusion coefficient comes from an ergodic Markov chain $(a_t^{\mathrm{e}})_{t\ge0}$ with invariant distribution $\pi$. %\swtd{ToDo: Motivation for this setting. For example, the diffusion coefficient may be unknown and information about it may be quantified through a posterior distribution from a (Bayesian) pre-stage experiment.}
To be more precise, %we assume that in general the diffusion coefficient is unknown, but we have access to some statistical information about $a$ in form of a probability distribution $\pi$. \sw{For example, this information may come from a pre-stage Bayesian experiment.} In many practical scenarios the explicit computation of $\pi$ or the generation of exact samples of $\pi$ is infeasible and 
we assume that this information about the diffusion coefficient has been generated by an MCMC algorithm. %This means, we assume that $(a_t^{\mathrm{e}})_{t\ge0}$ is an ergodic Markov chain with invariant distribution $\pi$. 
The observation model is then defined by 
\[u_t = S_t z^\ast + w_t\,, \]
where $S_t = \mathcal O \circ G_t$ with $G_t = G_{a_t^{\mathrm{e}}}$.

Let us now consider how we generate our correlated data. %In order to do so we can generate 
We consider %an ergodic Markov chain \sw{$(a_t^{\mathrm{e}})_{t\ge0}$} whose 
a stationary distribution given as Gaussian distribution. %, is the distribution of interest. 
%We will assume that our correlated data
%is normally distributed. 
In order to generate the Markov chain $(a_t^{\mathrm{e}})_{t\ge0}$ with invariant distribution $\pi$, we implement a MCMC method, in particular the Metropolis-Hastings MCMC (MH-MCMC) algorithm. 

\begin{figure}[h!]
\centering
\includegraphics[scale=0.3]{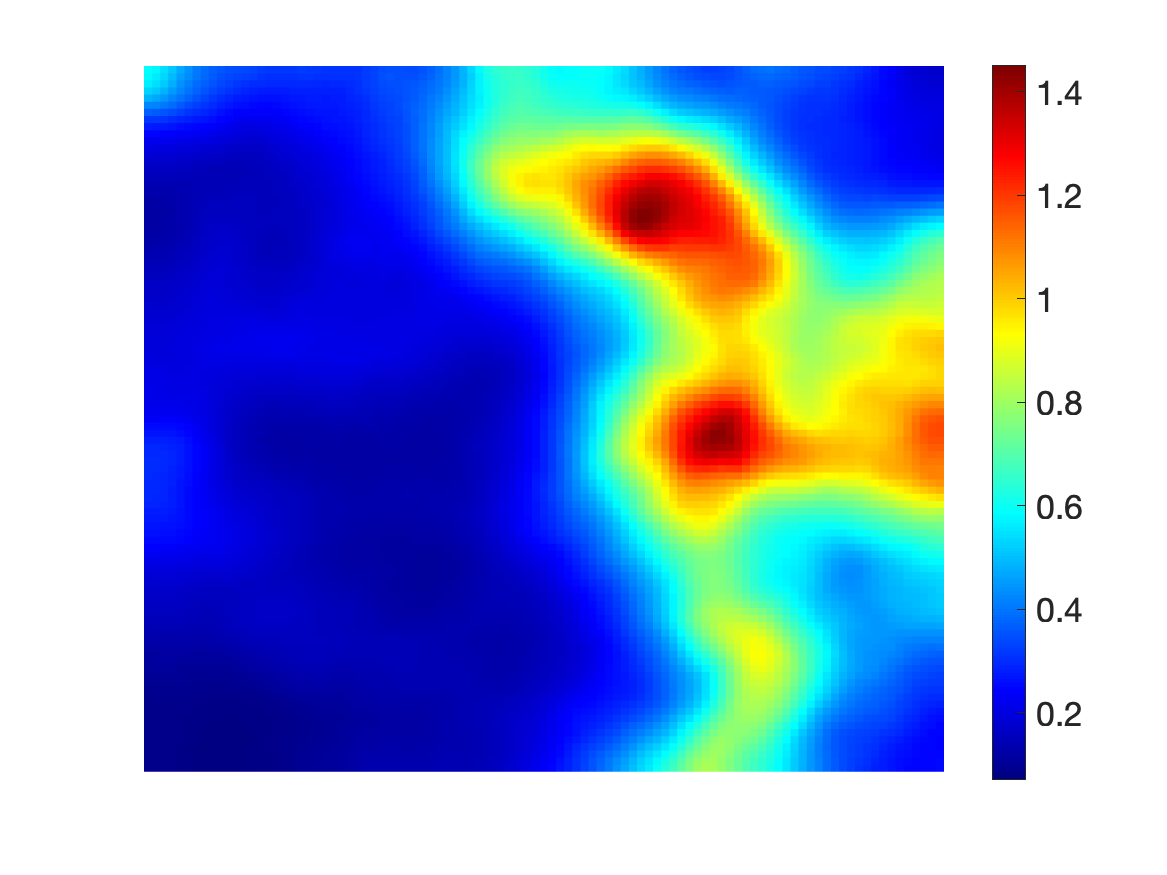}
\includegraphics[scale=0.3]{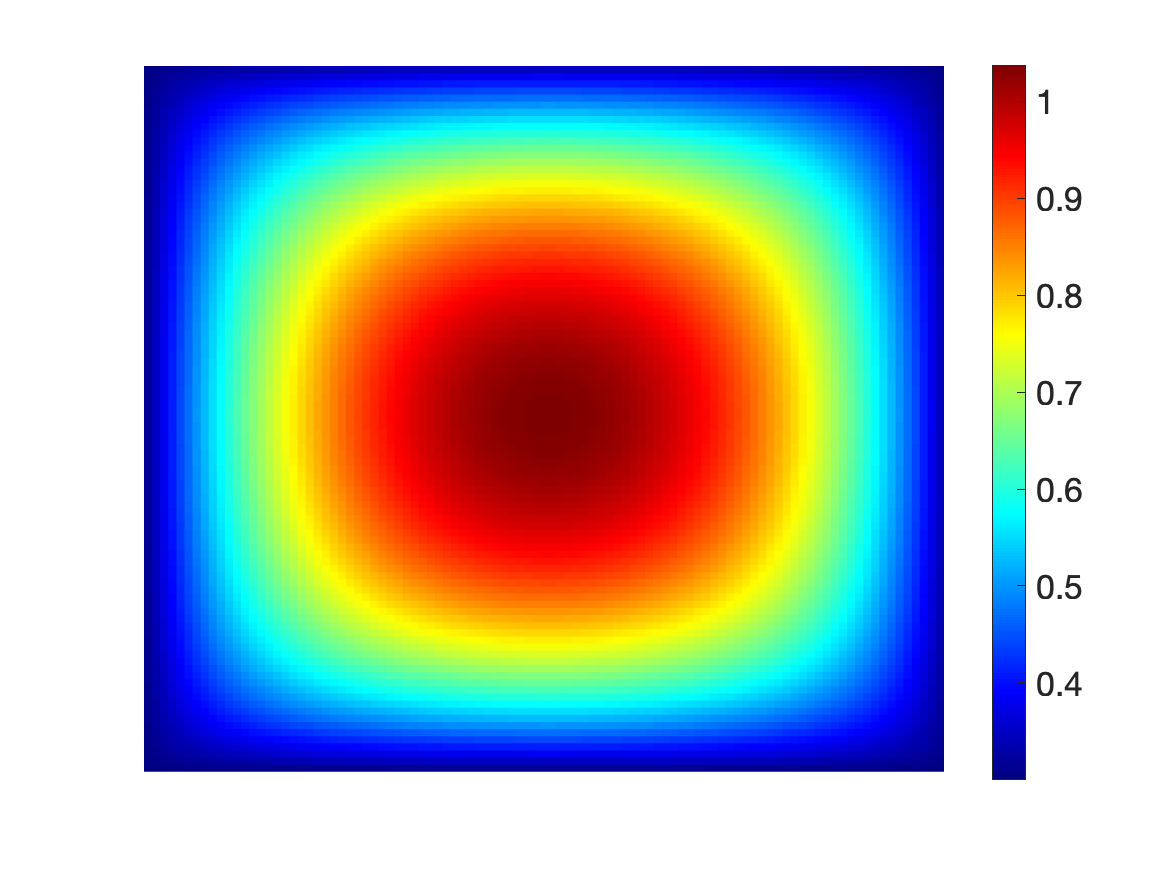}
\caption{ left:  Plot of our random field representation of $z^\ast$ based on the KLE \eqref{eq:kle}. right: Solution of the PDE \eqref{eq:darcy_flow}.}
\label{fig:fields}
\end{figure}

\begin{figure}[h!]
\centering
\includegraphics[scale=0.35]{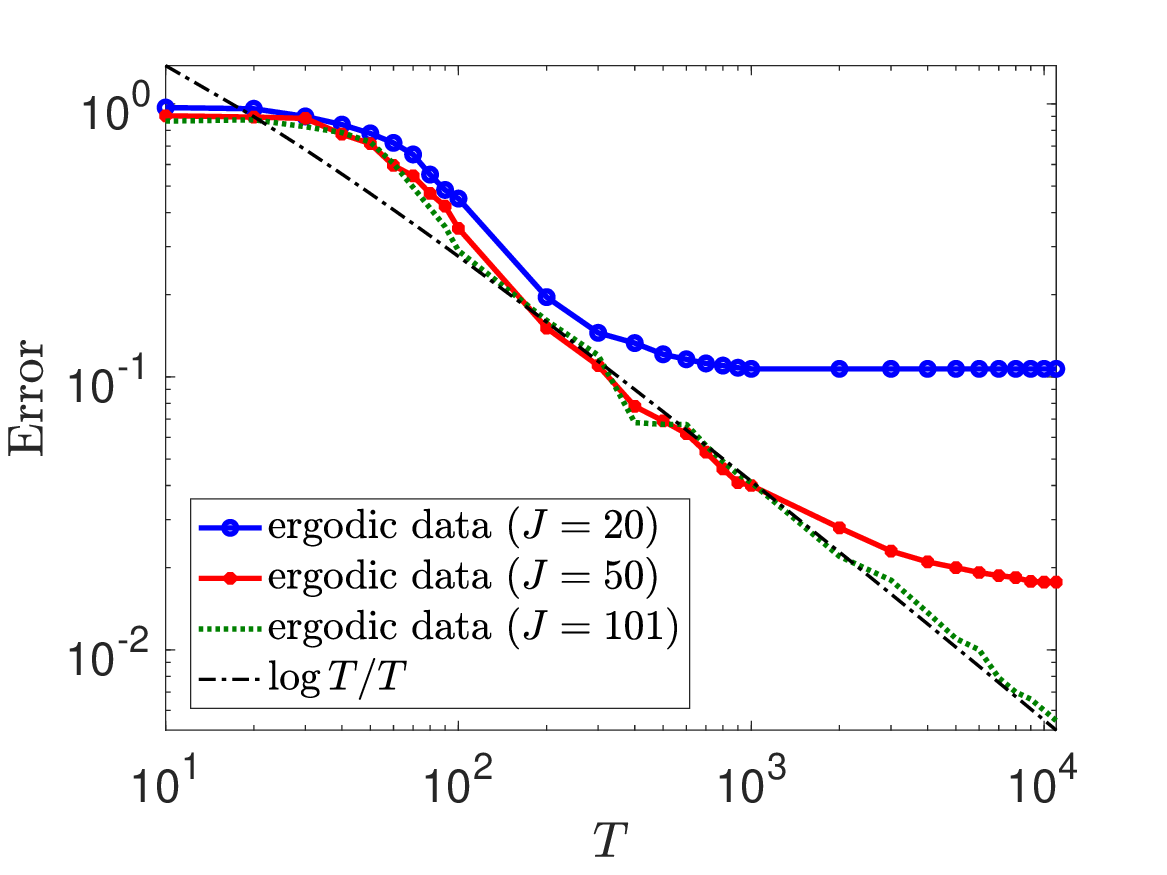}
\caption{Convergence plot for EKI using ergodic data. \textcolor{black}{We plot the error w.r.t.~the reference solution $z_{\mathrm{ref}}$ vs time $T$.}}
\label{fig:cov2}
\end{figure}

%We specify a Gaussian prior $\mu \sim \mathcal{N}(0,\mathcal{C})$ with covariance structure defined as a Mat\'{e}rn covariance function, which is simulated through the KLE, as defined in \eqref{eq:kle}. 
We use the MH-MCMC method based on proposing moves using on a pre-conditioned \textcolor{black}{Crank Nicolson} (pCN) scheme of the form, 
$$
y'=\sqrt{1-\beta^2}y+ \beta \epsilon, \quad \epsilon \sim \mathcal{N}(0,\mathcal{C}),
$$
with covariance structure defined as a Mat\'{e}rn covariance function, which is simulated through the KLE, as defined in \eqref{eq:kle}. The proposal distribution is then of the form 
$$
y'\sim \mathcal N((1-\beta^2)y,\beta^2 \mathcal C),
$$
 which uses an acceptance probability %based on whether 
 to either accepted or reject the proposed moves, defined as
%$$
%\alpha_{\textrm{acc}}(y,y') = \min\bigg\{1,\frac{\pi(y|y')\pi(y')}{\pi(y'|y)\pi(y)}\bigg\}\, .
%$$
$$
\alpha_{\textrm{acc}}(y,y') = \min\bigg\{1,\frac{\pi(y')}{\pi(y)}\bigg\}\, .
$$
%\swtd{What is the choice of $\pi$? Is it $\pi = \mathcal N(0,\mathcal C)$? And $\pi(y|y')$ are the proposal probabilities, right? pCN is symmetric, so it should hold $\pi(y|y') = \pi(y'|y)$, and we just need to consider $\pi(y')/\pi(y)$?  Yes Simon I have modified this, I just wanted to give the full expression}
For our numerical experiments we set $\beta=0.9$ of our proposal to ensure our acceptance rate is approximately $\alpha_{\mathrm{accep}} = 0.234$, which is consistent with the ``optimal" value for RWMH. When running our DEKI algorithm we again specify $J=50$ ensemble members where we place $T=10,000$. Our \textcolor{black}{ground truth} $z^{\ast}$ is again a realization of \eqref{eq:kle} with specific choices $(\ell,\nu) = (0.07,3.4)$. %\swtd{I have changed the description of $z^\ast$. Is this correct? Yes}
% based on a high-resolution MCMC simulation for our data. 
Figure \ref{fig:fields} presents our generated \textcolor{black}{ground truth} $z^{\ast}$ and the corresponding solution to the Darcy flow PDE. \textcolor{black}{Our step size for the experiment is chosen again as $\eta_t=t^{-1}$. } %\swtd{Is this choice $\eta_t$ left from the MCGD method? For EKI we don't need a step size, right? Just the scaling $h$.}%Furthermore for our pCN proposal, we set $\beta=0.9$. 
Our numerical experiments are shown in Figure \ref{fig:cov2} where we, as before, observe the theoretical rates %. Theoretically we observe the rates 
from Theorem~\ref{thm:main}.
\textcolor{black}{
Interestingly what we see is that at the beginning of the learning process, the ergodic data is slower as the initialization of the MCMC process induces a burn-in period. To allow for this we set the initialization as normal distribution with mean $\mu_{\textrm{mean}} = 3$. We again compare the effect of the ensemble size. For the cases of $J=\{20,50\}$, i.e. when $J<d$ we do not observe convergence, similar to the periodic and i.i.d. data. However, when we specify $J=101>100=d$ we see no plateau, implying convergence  Theorem~\ref{thm:main}. In order for us to ensure convergence, for all types of data considered, we would require techniques that would break the subspace property such as localization.}
 
\begin{figure}[h!]
\centering
\includegraphics[scale=0.35]{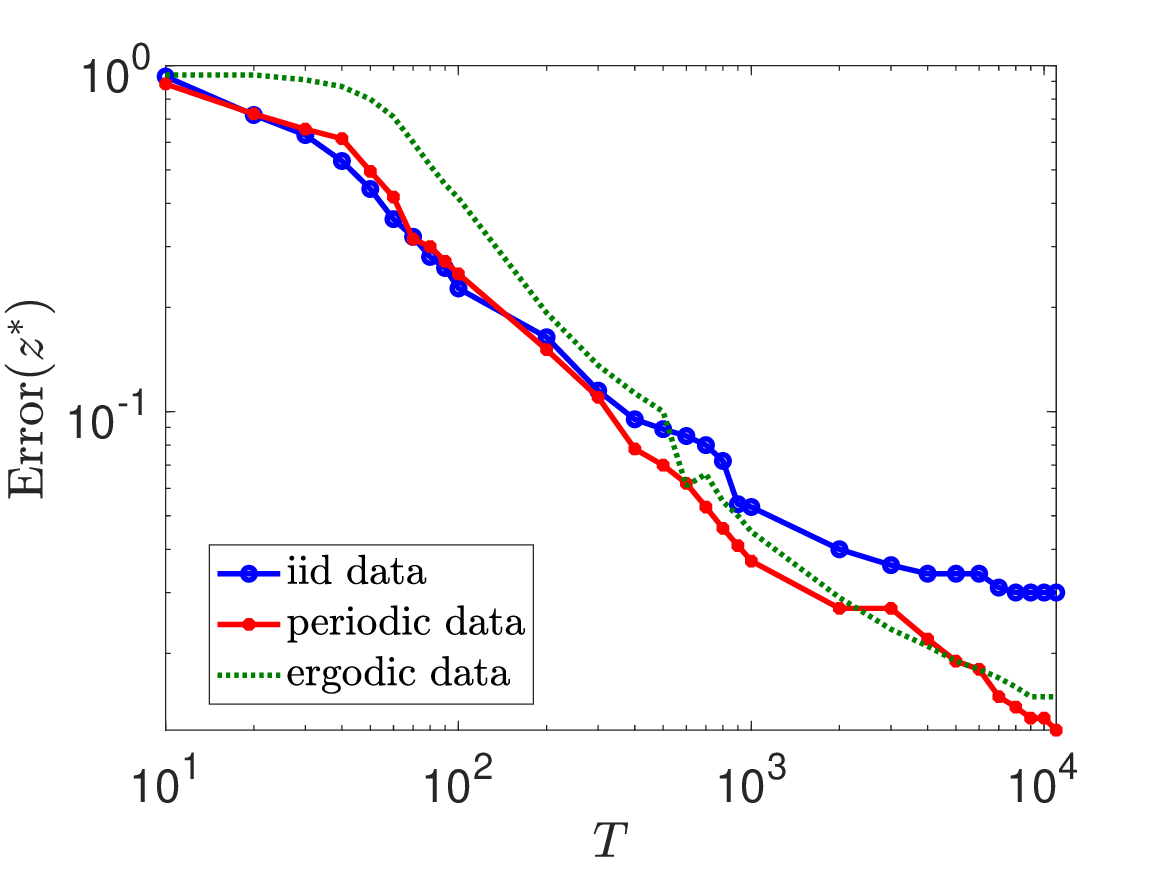}
\caption{\textcolor{black}{Convergence plot for EKI using  i.i.d., periodic and ergodic data. \textcolor{black}{We plot the error w.r.t.~the ground truth $z^\ast$ vs time $T$.}}}
\label{fig:revision}
\end{figure}

\paragraph{\bf Ground truth error}
\textcolor{black}{As a final experiment is to consider the same convergence plots, as before, where now the 
 error is w.r.t. the ground truth, i.e. $\E[\|\bar z_t - z^{\ast}\|]$. This differs to the theoretical results obtained in Section \ref{sec:EKI}, where our setup is similar to the one described in Section \ref{sec:num}. Our simulations are presented in Figure \ref{fig:revision}
 which demonstrate that the error does decrease as we increase $T$. This is the case for both periodic and ergodic data. Further work is required on proving convergence to the ground truth based on regularization theory.}
 
\section{Conclusion}
\label{sec:conc}

Within the field of inverse problems, it is common to learn unknown parameters 
from noisy data where the forward operator is static. In this work we consider the setting 
where instead we aim to learn parameters where there is a time-dependent forward operator, i.e.
the forward model is dynamic. We consider this in the setting of an optimizer for black-box inverse problems, which is the EKI.
Our new algorithm, which we refer to as \textcolor{black}{DEKI}, is introduced in a modified setting where we assume a linear least squares problem with the
addition of Tikhonov regularization.  A number of important results are derived which include a moment bound
and ensemble collapse result, and our main result which is a convergence analysis. Our results apply to
different data cases such as (i) \textit{i.i.d. data}, (ii) \textit{ergodic data} and (iii) \textit{periodic data}. 
Numerical experiments are conducted on a toy model %, a toy matrix problem, which is similar to a block co-ordinate minimization problem 
based on a linear elliptic PDE-constrained optimization problem using the 2D Darcy flow model. Our experiments demonstrate and verify our theoretical findings. 

We conclude this article with a number of potentially useful directions to enhance our work. 
These are summarized below.
\begin{itemize}
\item Arguably the most important extension is related to theory, in terms of two aspects. The first being a 
non-linear analysis \cite{CT22,W2022}, and the other being the extension to an infinite-dimensional analysis. Both works, in particular the first, are
of significant challenge as there has been limited work, in general, for the linear case which we are able
to provide analysis. The latter work is currently ongoing work by the authors.
\item An obvious question is if this work could be extended to adaptive choice of $\alpha_t$? For 
example one could consider bilevel learning \cite{CSTW22}, as well as other adaptive Tikhonov strategies \cite{WCS22}.
%\item Exploit more advanced methodologies, such as momentum based optimization,
%or methods that use Hessian information, such as Newton solvers? This is ongoing work. 
\item Another interesting direction to consider is exploiting
this framework in the continuous-time version. In the context of continuous dynamics with ergodic data, one could perhaps
using Langevin based ideas from MALA and HMC, as examples for advanced proposals within MH methods.
\item A final consideration, \textcolor{black}{related to Remark \ref{rem:ensemble}}, could be on the influence of either variance inflation or localization \cite{NKC18,CCS21,TMK16,TM2023}.
{Alternatively, the dropout technique from machine learning community can also be adapted on EKI for $J<d$ scenarios \cite{LRT23}.} In the context of \textcolor{black}{DEKI}, it may
also be the case where one could attain improved rates of convergence.  {At the same time, the numerical experiments in Section \ref{sec:num} use $J<d$ ensembles without the aforementioned modifications. Somehow their performances are satisfactory. This may be because the effective dimension of these problems are low \cite{GS24,MT17} and DEKI can implicitly exploit it. However, further studies are required to understand this phenomena.}
%\item Thus far we implemented a time-dependent version of EKI, which to the best of our knowledge is the first
%time-dependent derivative-free optimizer. A natural question is how does this compare to other DFO's. One extension
%that could be done is to use this work, prehaps for particle swarm optimization or the CMA-ES method.
%\swtd{as stated above, in the linear setting we can not sell EKI as derivative-free optimizer, since we are computing the exact gradients through $S_t^\top$.}
\end{itemize}

\section*{Acknowledgments}
NKC is supported by an EPSRC-UKRI AI for Net Zero Grant: “Enabling CO2
Capture And Storage Projects Using AI”, (Grant EP/Y006143/1).
The work of XTT has been funded by Singapore MOE grant A-8000459-00-00. {The authors are very grateful for helpful discussions with Claudia Schillings.}


\begin{thebibliography}{1}

\bibitem{AK23}
H. Albers and T. Kluth.
Time-dependent parameter identification in a Fokker-Planck equation based magnetization model of large ensembles of nanoparticles.
Arxiv preprint, arxiv:2307.03560, 2023.


%\bibitem{AD13}
%A. Agarwal and J. C. Duchi.
%The Generalization Ability of Online Algorithms for dependent data.
%\newblock{\em IEEE Transactions on Information Theory}, 59(1), 2013.


%\bibitem{BRS18}
%J. Bhandari, D. Russo, and R/ Singal. 
%A finite time analysis of temporal difference learning with linear function approximation. 
%\newblock{\em Conference on Learning Theory}, 1691--1692, PMLR, 2018.
\bibitem{BB18}
M. Benning and M. Burger.
\newblock{Modern regularization methods for inverse problems}, {\em Acta Numerica}, {27},
\newblock 2018.

\bibitem{BSW18}
D. Bl\"{o}mker, C. Schillings, P. Wacker.
 \newblock{A strongly convergent numerical scheme from ensemble Kalman inversion,}
 \newblock{\em SIAM J. Numerical Analysis}, 56(4), 2018.


%\bibitem{BHG03}
%L. T. Biegler, M. Heinkenschloss, O. Ghattas and B. B. Waanders.
%Large-scale PDE-constrained optimization.
%\newblock{\em Lecture Notes in Computational Science and Engineering}, Springer, Berlin, 2003.



\bibitem{BSWW19}
D. Bl\"{o}mker, C. Schillings, P. Wacker and S. Weissmann.
\newblock{Well posedness and convergence analysis of the ensemble Kalman inversion},
\newblock{\em Inverse Problems}, 2019.

\bibitem{BSWW22}
D. Bl\"{o}mker, C. Schillings, P. Wacker and S. Weissmann.
\newblock{Continuous Time Limit of the Stochastic Ensemble Kalman Inversion: Strong Convergence Analysis},
\newblock{\em SIAM Journal on Numerical Analysis}, 60(6), 2019.


\bibitem{NKC18}
N.~K. Chada.
\newblock Analysis of hierarchical ensemble Kalman inversion.
\newblock{arXiv preprint arXiv:1801.00847}, 2018.

\bibitem{CCS21}
N. K. Chada, Y. Chen, D. Sanz-Alonso.
Iterative ensemble Kalman methods: A unified perspective with some new variants.
\newblock{\em Foundations of Data Science}, 3(3), 331--369, 2021.



\bibitem{CIRS18}
N.~K Chada, M.~A. Iglesias, L. Roininen, and A.~M. Stuart.
Parameterizations for ensemble {K}alman inversion.
\newblock {\em Inverse Problems}, 34(5):055009, 2018.



\bibitem{CSW19}
N.~K. Chada, C. Schillings, and S. Weissmann.
On the incorporation of box-constraints for ensemble {K}alman inversion.
\newblock {\em Foundations of Data Science}, 1(2639-8001 2019 4 433):433--456,
  2019.


\bibitem{CST19}
N.~K. Chada, A.~M. Stuart and X.~T. Tong.
\newblock Tikhonov regularization within ensemble {K}alman inversion.
\newblock {\em SIAM Journal on Numerical Analysis}, 58(2):1263--1294, 2020.



\bibitem{CT22}
N.~K. Chada and X.~T. Tong.
\newblock Convergence acceleration of ensemble {K}alman inversion in nonlinear
  settings.
\newblock{\em Math. of Comp.}, 91(335), 1247--1280, 2022.



\bibitem{CSTW22}
N. K. Chada, C. Schillings, X. T. Tong and S. Weissmann.
Consistency analysis of bilevel data-driven learning in inverse problems.
\newblock{\em Communications in Mathematical Sciences},  20(1), 123--164, 2022.


\bibitem{CMT14}
N. Chen, A. J. Majda and X. T. Tong.
Information barriers for noisy Lagrangian tracers in filtering random incompressible flows.
\newblock{\em Nonlinearity}, 27,2133--2163, 2014.


%\bibitem{CR21}
%P. Chen and J.  O. Royset.
%Performance bounds for PDE-Constrained optimization under uncertainty.
%Arxiv preprint, arxiv:2110.10269, 2021.

%\bibitem{CE17}
%J. Chung and M. I Espanol.
%\newblock{Learning regularization parameters for general-form Tikhonov}.
%\newblock{\em Inverse Problems}, Volume 33, Number 7, 2017.



\bibitem{DL21}
Z. Ding and Q. Li.
\newblock{Ensemble Kalman sampler: mean-field limit and convergence analysis},
\newblock{\em SIAM J. Math. Anal.}, 53(2), 1546--1578, 2021.


%\bibitem{DT22}
%J.~Dong and X.~T. Tong. 
%\newblock{Stochastic Gradient Descent with Dependent Data for Offline Reinforcement Learning}.
%\newblock{\em arXiv Preprint}, 2022.

%\bibitem{CLTZ16}
%X.Chen, J. D. Lee, X. T. Tong and Y. Zhang.
%\newblock{Statistical inference for model parameters in stochastic gradient descent}.
%\newblock{arXiv preprint arXiv:1610.08637}, 2016.



%\bibitem{DAJ12}
%J. C. Duchi, A. Agarwal, M. Johansson, and M. I. Jordan.
%Ergodic mirror descent.
%\newblock{\em SIAM J. Optim.}, 22(4), 2012.


%\bibitem{DJM20}
%A. Durmus, P. Jiménez, E. Moulines, S. Said and Hoi-To Wai.
%Convergence analysis of Riemannian stochastic approximation schemes.
%\newblock{\em arXiv Preprint, arxiv:005.13284}, 2020.


%
% \bibitem{EHN96}
% H.W. Engl, K. Hanke and A. Neubauer.
% \newblock{Regularization of inverse problems}, {\em Mathematics and its Applications}, Volume {375}, Kluwer Academic Publishers Group, Dordrecht,
% \newblock 1996.
 

\bibitem{GE09}
G. Evensen.
\newblock {\em Data Assimilation: The Ensemble Kalman Filter}. 
\newblock{Springer}, 2009.


\bibitem{GE03}
G. Evensen. 
\newblock{The ensemble Kalman filter: Theoretical formulation and practical implementation.}
\newblock{\em Ocean dynamics}, 53(4):343--367, 2003


\bibitem{ME23}
M. Even.
Stochastic gradient descent under Markovian sampling schemes.
\newblock{\em Proceedings of the 40th International Conference on Machine Learning}, 
PMLR 202:9412--9439, 2023.



\bibitem{GHLS19}
A. Garbuno-Inigo, F. Hoffmann, W. Li and A. M. Stuart,
\newblock{Gradient structure of the ensemble Kalman flow with noise.}
\newblock{SIAM J. Applied Dynamical Systems}, 19(1), 412--441, 2020.
\textcolor{black}{
\bibitem{GS24}
O. Al-Ghattas and D. Sanz-Alonso
Non-asymptotic analysis of ensemble Kalman updates: effective dimension and localization.
\newblock{\em Information and Inference: A Journal of the IMA}, 13(1), 2024.}




\bibitem{HLS23}
M. Hanu, J. Latz and C. Schillings.
\newblock{Subsampling in ensemble Kalman inversion.}
\newblock{\em Inverse Problems}, 39(9) 2023.
\textcolor{black}{
\bibitem{HOS21}
A. Hauptmann, O. Öktem, and C. Schönlieb.
Image reconstruction in dynamic inverse problems with temporal models.
K. Chen et al. (eds.), \newblock{\em Handbook of Mathematical Models and Algorithms in Computer
Vision and Imaging}, 2021.}

\bibitem{MAI16}
M. A. Iglesias.
\newblock {A regularising iterative ensemble Kalman method for PDE-constrained inverse problems.} 
{\em Inverse Problems}, {32}, 2016.


\bibitem{ILS13}
M. A. Iglesias, K. J. H. Law and A. M. Stuart.
\newblock {ensemble Kalman methods for inverse problems}. 
\newblock {\em Inverse Problems}, {29}, 2013.


\bibitem{IY21}
M. A. Iglesias and Y. Yang.
Adaptive regularisation for ensemble Kalman inversion.
\newblock{\em Inverse Problems}, 37 025008, 2021.

\bibitem{IPT19}
M. Iglesias, M. Park and M. V. Tretyakov.
Bayesian inversion in resin transfer molding.
\newblock{\em Inverse Problems}, 34(10), 105002, 2019.


%\bibitem{KS04}
%J. Kaipio and E. Somersalo.
%\newblock {\em Statistical and Computational Inverse problems.} Springer Verlag, New
%York,
%\newblock 2004.


\bibitem{BK17}
B. Kaltenbacher.
All-at-once versus reduced iterative methods for time dependent inverse problems.
\newblock{\em Inverse Problems}, 33, p. 064002, 2017.


\bibitem{KSW21}
B. Kaltenbacher, T. Schuster, and A. Wald.
\newblock{\em Time-dependent Problems in Imaging and Parameter Identification}.
Springer International Publishing, Cham, 2021.


\bibitem{KSW21b}
R. Klein, T. Schuster, and A. Wald.
Sequential subspace optimization for recovering stored energy functions in hyperelastic materials from time-dependent data.
\newblock{\em In: Time-dependent Problems in Imaging and Parameter Identification}, 2021.
\textcolor{black}{
\bibitem{KMW15}
Y. Kwong, A. O. Mel, G. Wheeler and J. M. Troupis.
Four-dimensional computed tomography (4DCT): a review of the current status and applications. 
\newblock{\em J. Med. Imag. Radiat. Oncol.}, 59(5), 545--554, 2015.}
\bibitem{LSZ15}
K. J. H. Law, A. M. Stuart and K. Zygalakis.
\newblock{\em Data Assimilation: A Mathematical Introduction}.
\newblock{Texts in Applied Mathematics, Springer}, 2015. 
\textcolor{black}{
\bibitem{LS17}
A. Lechleiter and J. W. Schlasche. 
Identifying Lame parameters from time-dependent elastic wave measurements. 
\newblock{\em Inverse Problems in Science and Engineering}, 25, 2--26, 2017.
\bibitem{TR13}
T. Roubicek.
\newblock{\em Nonlinear Partial Differential Equations with Applications}.
International Series of Numerical Mathematics, Springer Basel, 2013.
}
\bibitem{LR09} 
G. Li  and A. C. Reynolds.
\newblock{Iterative ensemble Kalman filters for data assimilation}.
 \emph{SPE J} 14  496-505, 2009
{\bibitem{LRT23}
S. Liu, S. Reich and X. T. Tong. \newblock{Dropout Ensemble Kalman inversion for high dimensional inverse problems.} arXiv preprint arXiv:2308.16784. 2023}
\textcolor{black}{
 \bibitem{MT17} 
A. Majda and X. T. Tong. 
Performance of Ensemble Kalman Filters in Large Dimensions.
\newblock{\em Communications on Pure and Applied Mathematics}, 17(12), 892--937, 2017.
}
 \bibitem{MW06} 
A. Majda and X. Wang. 
\newblock{\em Non-linear Dynamics and Statistical Theories for Basic Geophysical Flows,} Cambridge University Press, 
\newblock 2006.

\bibitem{MT93}
S. P. Meyn  and R. L. Tweedie.
\newblock{\em Markov Chains and Stochastic Stability}.
Cambridge University Press, 1993.
%\bibitem{LLG15}
%B. Liu, J. Liu, M. Ghavamzadeh, S. Mahadevan and M. Petrik. 
%Finite sample analysis of proximal gradient TD algorithms. 
%\newblock{31st Conference on Uncertainty in Artificial Intelligence (UAI)}, 504--513, 2015.



%\bibitem{MN21}
%M.  Martin and F Nobile.
%PDE-constrained optimal control problems with uncertain parameters using SAGA.
%\newblock{SIAM/ASA J. Uncertain. Quantif.}, 9(3), 979--1012, 2021.


\bibitem{TTN19}
T.T.N. Nguyen.
Landweber–Kaczmarz for parameter identification in time-dependent inverse problems: all-at-once versus reduced version.
\newblock{\em Inverse Problems}, 35, 035009, 2019.
\textcolor{black}{
\bibitem{ORL08}
D. Oliver, A. Reynolds and N. Liu. 
\newblock{\em Inverse Theory for Petroleum Reservoir Characterization and History Matching}.
Cambridge University Press, 2008.}

\bibitem{SCI92} 
E. Somersalo, M. Cheney and D. Isaacson.
\newblock{Existence and Uniqueness for Electrode Models for Electric Current Computed Tomography}, {\em SIAM J. Appl. Math.,} {52}, 1023-1040, 
\newblock 1992.
%\bibitem{SR19}
%S. Reich.
%\newblock{ Data Assimilation}.
%\newblock{\em Acta Numerica, in preparation}, 2019.



\bibitem{SS17}
C.~Schillings and A.~M. Stuart.
\newblock{Analysis of the ensemble {K}alman filter for inverse problems}.
\newblock {\em SIAM J. Numer. Anal.}, 55(3):1264--1290, 2017.


%\bibitem{VHS14}
%V. H. Schulz. 
%A Riemannian view on shape optimization. 
%\newblock{\em Foundations of Computational Mathematics}, 14, 483--501, 2014.


%\bibitem{SSW15}
%V. H. Schulz, M. Siebenborn and K. Welker. 
%Towards a Lagrange-Newton approach for PDE constrained shape optimization. In A. Pratelli and G. Leugering, editors, 
%\newblock{\em New Trends in Shape Optimization, volume 166 of International Series of Numerical Mathematics},
%229--249. Springer, 2015.


\bibitem{AMS10}
A. M. Stuart.
\newblock{Inverse problems: A Bayesian perspective}. 
\newblock{\em Acta Numerica}, Vol. {{19}}, 451-559, 2010.


\bibitem{SL19}
T. Sun and D. Li.
Decentralized Markov chain gradient descent.
arxiv preprint, arXiv:1909.10238, 2019.



\bibitem{SSY18}
T. Sun, Y. Sun and W. Yin.
On Markov chain gradient descent.
\newblock{\em 32nd Conference on Neural Information Processing Systems}, 2018.


\bibitem{AT87}
A. Tarantola. 
\newblock{\em Inverse Problem Theory and Methods for Model Parameter Estimation}. 
Elsevier, 1987.


%\bibitem{TKX12}
%H. Tiesler, R. M. Kirby, D. Xiu, and T. Preusser.
%Stochastic collocation for optimal control problems with stochastic PDE constraints. 
%\newblock{\em SIAM J. Control Optim.}, 50, 2659--2682, 2012.



\bibitem{TMK16} 
X. T. Tong, A. J. Majda and D. Kelly.  
\newblock{Nonlinear stability of the ensemble Kalman filter with adaptive covariance inflation.}
\newblock{\em Commun. Math. Sci.}, 14(5):1283--1313, 2016.

\bibitem{TM2023}
X. T. Tong and M. Morzfeld
\newblock{Localized ensemble Kalman inversion}
\newblock{\em Inverse Problems}, 32(1) 064002, 2023.


%\bibitem{FT10}
%F. Tröltzsch.
%\newblock{\em Optimal Control of Partial Differential Equations: Theory, Methods, and Applications}. 
%Graduate Studies in Mathematics, American Mathematical Society, 2010.


\bibitem{LH16}
T. van Leeuwen and F. J. Herrmann.
A penalty method for PDE-constrained optimization in inverse problems.
\newblock{\em Inverse Problems}, 39(6) 015007, 2016.


\bibitem{WLY23}
P. Wang, Y. Lei, Y. Ying and D-X. Zhou.
\newblock{Stability and generalization for Markov Chain stochastic gradient methods. }
In \newblock{\em Advances in Neural Information Processing Systems}, 2022.


\bibitem{WCS22}
S. Weissmann, N. K. Chada, C. Schillings and X. T. Tong.
\newblock{Adaptive Tikhonov strategies for ensemble Kalman inversion.}
\newblock{\em Inverse Problems}, 38(4), 2022.


\bibitem{majda2015intermittency}
A. Majda and X.T. Tong,
\newblock{Intermittency in turbulent diffusion models with a mean gradient}
\newblock{\em Nonlinearity}, 28(11), 2015.



\bibitem{W2022}
S. Weissmann
\newblock{Gradient flow structure and convergence analysis of the ensemble Kalman inversion for nonlinear forward models}
\newblock{\em Inverse Problems}, 38(10), 2022.


\end{thebibliography}
\end{document}